\documentclass{article}
\usepackage{amssymb}
\usepackage{amsmath}
\usepackage{latexsym}

\newlength{\mytopmargin}
\newlength{\myleftmargin}
\def\rlx{\relax\leavevmode}
\def\inbar{\vrule height1.5ex width.4pt depth0pt}
\def\zz{\rlx\hbox{\small \sf Z\kern-.4em Z}}
\def\rr{\rlx\hbox{\scriptsize \rm I\kern-.18em R}}
\def\nn{\rlx\hbox{\rm I\kern-.18em N}}
\def\qq{\rlx\hbox{\,$\inbar\kern-.3em{\rm Q}$}}

\newcommand{\ka}{\kappa}

\newcommand{\de}{\delta}
\newcommand{\om}{\omega}

\newtheorem{lemma}{Lemma}[section]

\newtheorem{cor}[lemma]{Corollary}
\newtheorem{prop}[lemma]{Proposition}

\begin{document}

\noindent
\begin{center}{  \Large\bf
Symmetric and Non-symmetric Macdonald Polynomials }
\end{center}
\vspace{5mm}
 
\noindent
\begin{center} Dan Marshall\footnote{email: danm@maths.mu.oz.au; supported
by an APA Scholarship.} \\[2mm]
{\it Department of Mathematics and Statistics, \\
 University of Melbourne,\\
 Parkville, Victoria 3052, Australia}
\end{center}
\vspace{.5cm}
		 
\date{}

\begin{quote}
The symmetric Macdonald polynomials are able to be constructed
out of the non-symmetric Macdonald polynomials. This allows us
to develop the theory of the symmetric Macdonald polynomials
by first developing the theory of their non-symmetric 
counterparts. In taking this approach we are able to obtain
new results as well as simpler and more accessible derivations
of some of the known fundamental properties of both kinds of polynomials.
\end{quote}

\section{Introduction}

The symmetric Macdonald polynomial $P_{\kappa}:=P_{\kappa}(x;q,t)$
is a polynomial of $n$ variables \(x=(x_{1}, \ldots ,x_{n})\) having 
coefficients
in the field $\mathbb{Q}(q,t)$ of rational functions in $q$ and $t$. The symbols
$q$ and $t$ can be interpretated either as indeterminants
or as parameters ranging over
$0<q,\, t<1$. 
The symmetric Macdonald polynomial $P_{\ka}(x;q,t)$ is labeled by a
partition of length $\leq n$ and can be defined as the unique 
eigenfunction of the operator
\begin{equation} \label{1.1}
D_{n}^{1}(q,t)=\sum_{i=1}^{n}\sum_{i\neq j}
\frac{tx_{i}-x_{j}}{x_{i}-x_{j}}\tau_{i}
\end{equation}
which is of the form
\begin{equation} \label{1.2}
P_{\kappa}(x;q,t)=m_{\kappa}(x)+
\sum_{\mu<\kappa}u_{\kappa\mu}m_{\mu}(x)
\end{equation}
In (\ref{1.2}), $m_{\kappa}(x)$ is the monomial symmetric 
function in variables
 $x_{1},\ldots ,x_{n}$ and the sum is over the partitions $\mu $ 
which have the same
modulus as $\kappa $, but are smaller in dominance ordering. 
The $q$-shift 
operator $\tau_{i}$ in (\ref{1.1}) acts on 
functions so that
\begin{equation}
(\tau_{i}f)(x_{1}\ldots x_{n})=f(x_{1},\ldots ,qx_{i},\ldots x_{n})
\end{equation}
The symmetric Macdonald polynomials have been the subject of 
much recent study, both for their mathematical properties 
\cite{BF97g},\cite{KS96},\cite{Ok97}
and their applications to the the trigonometric Ruijsenaars-Schneider
quantum many body model \cite{Ko96}. They can be viewed as a 
$q$-generalisation of the symmetric Jack polynomials, 
the latter being obtained from the former in the limit 
$q\rightarrow1$ with 
$q=t^{\alpha}$ and $\alpha$ fixed. In this paper we will 
develop the theory of the Macdonald polynomials by 
generalising the approach taken by Baker and Forrester 
\cite{BF97c} towards the Jack polynomials.

The strategy is to first develop the theory of non-symmetric
Macdonald polynomials. These polynomials were first 
introduced \cite{Ch95a},\cite{Ma95a} some 
time after the seminal work of 
Macdonald \cite{Ma95}
on the symmetric Macdonald polynomials. The symmetric 
polynomials can be constructed from their non-symmetric 
counterparts. This opens the way to using the theory 
of the non-symmetric Macdonald polynomials to develop 
the theory of the symmetric Macdonald polynomials. 
In taking this approach we will obtain new results 
as well as new and simpler derivations of known
results. In the later case references will be 
provided to the original contributors.

\setcounter{equation}{0}
\section{Preliminaries}

In this section we will revise the basic definitions 
and results of the non-symmetric, symmetric and 
$q$-antisymmetric polynomials.
Following Macdonald \cite{Ma95a}, the symmetric and 
$q$-antisymmetric
polynomials will be constructed in terms of 
their non-symmetric counterparts, rather 
than as an independent entity as would stem 
from making (\ref{1.1}) and (\ref{1.2}) the starting point.
In addition, dual non-symmetric Macdonald polynomials
will be defined and related to the symmetric 
and $q$-antisymmetric Macdonald polynomials.
The results presented on this topic are for the
most part new.

The non-symmetric Macdonald polynomials are defined in terms of
operators which generate an extended affine Hecke algebra 
(see e.g.\ \cite{KN}). Let $s_{ij}$ be the operator which acts on functions of
\mbox{$x:=(x_1,\ldots, x_n)$} by interchanging
the variables $x_{i}$ and $x_{j}$. 
The Demazure-Lustig operators are defined by
\begin{eqnarray} \label{pucallpa}
T_i &:=& t + \frac{tx_i-x_{i+1}}{x_i - x_{i+1}}\left(s_i -1 \right)
\hspace{0.5cm} i=1,\ldots,n-1 \hspace{0.5cm} \mbox{and}\\
T_0 &:=& t + \frac{qtx_n-x_1}{qx_n - x_1}\left(s_0 -1 \right)
\end{eqnarray}
where $s_{i}:= s_{i\, i+1}$ and $s_{0} := s_{1n}\tau_1\tau_n^{-1}$.
The operators  $T_i$ have the following action on the monomial
 $x_i^a\,x_{i+1}^b$
for $1\leq i\leq n-1$ (see e.g \cite{KN}).
\begin{equation} \label{t.action}
T_i\, x_i^a x_{i+1}^b = \left\{ \begin{array}{ll}
(1-t)x_i^{a-1}x_{i+1}^{b+1} +
\cdots +(1-t)x_i^{b+1}x_{i+1}^{a-1} + x_i^bx_{i+1}^a& a > b \\
t x_i^ax_{i+1}^a & a=b \\
(t-1)x_i^{a}x_{i+1}^{b} +
\cdots +(t-1)x_i^{b-1}x_{i+1}^{a+1} + t x_i^bx_{i+1}^a & a < b
\end{array} \right.
\end{equation}
The operator $\om $ is defined by
\begin{equation} 
\om := s_{n-1}\cdots s_2\,s_1\tau_1 =
s_{n-1}\cdots s_i\tau_i s_{i-1}\cdots s_1 
\end{equation}
The extended affine Hecke algebra is then generated by elements $T_i$,
$0\leq i\leq n-1$ and $\om$, satisfying the relations
\begin{eqnarray}
(T_i-t)\,(T_i+1) &=& 0 \label{tdefs.1}\\
T_i\;T_{i+1}\;T_i &=& T_{i+1}\;T_i\;T_{i+1} \label{tdefs.2}\\
T_i\;T_j &=& T_j\;T_i \qquad |i-j| \geq 2 \\
\om\;T_i &=& T_{i-1}\;\om  \label{tdefs.4}
\end{eqnarray}
where the indices $0,1,\ldots,n-1$ are understood as 
elements of ${\mathbb{Z}}_{n}$.
From the quadratic relation (\ref{tdefs.1}), we have the identity
\begin{equation} 
T_{i}^{-1}=t^{-1}-1+t^{-1}T_{i}
\end{equation}
Given a permuation $\sigma$ with reduced word decomposition
$\sigma :=s_{i_{1}}\cdots s_{i_{p}}$
we define
\begin{equation}
T_{\sigma}:=T_{i_{1}}\cdots T_{i_{p}}
\end{equation}
The Cherednik operators \cite{ch91},\cite{Ch94} are defined by
\begin{equation}
Y_i = t^{-n+i}\,T_i\cdots T_{n-1}\;\om\;T_1^{-1}\cdots T_{i-1}^{-1},
\hspace{2cm} 1\leq i\leq n
\end{equation}
The fact that the Cherednik operators commute with 
each other, along with the triangularity of their action on 
$x^{\eta}:=x_{\eta_{1}} \cdots x_{\eta_{n}}$,
implies that they possess a set of simultaneous eigenfunctions. 
These are the
non-symmetric Macdonald polynomials $E_{\eta}$ which can be 
defined by the conditions
\begin{eqnarray} \label{E1}
E_{\eta}(x;q,t) &=& x^{\eta} + 
\sum_{\nu\prec \eta} b_{\eta\nu} x^{\nu} \\ \label{E2}
Y_i \, E_{\eta}(x;q,t) &=& t^{\bar{\eta}_i}\,E_{\eta}(x;q,t) \qquad
1\leq i \leq n
\end{eqnarray}
where
\begin{equation}\label{e-val.1}
\bar{\eta}_i := {\lambda^{-1}}\eta_i - l^{'}_{\eta}(i) \qquad 
l^{'}_{\eta}(i)
:=\#\{k<i\;|\;\eta_k \geq \eta_i\} +\#\{k>i\;|\;\eta_k > \eta_i\}
\end{equation}
with $\lambda$ is  parameter such that $t=q^{\lambda} $.
Let  $\eta^{+}$ be the unique partition obtained by permuting
$\eta$.
The partial order $\prec$ is defined on compositions having the same
modulus so that
\begin{equation}
\nu\prec\eta \quad\mbox{iff}\quad \nu^+ < \eta^+
\quad\mbox{or in the case $\nu^+ = \eta^+$}\quad \nu < \eta
\end{equation}
where $<$ is the usual dominance ordering for $n$-tuples, that is,
 $\nu < \eta$ iff $\sum_{i=1}^p (\eta_i - \nu_i) \geq 0$,
for all $1\leq p \leq n$.

Following Sahi \cite{Sa96} $l^{'}(s):=l^{'}(i)$ is called the 
leg colength of the node
$s=(i,j)$ in the composition $\eta$. The arm length 
$a(s) $, arm colength ${a'}(s)$
and leg length $l(s)$ are defined by
\begin{eqnarray}
a(s)= \eta_i - j && l(s) = 
\#\{k>i|j\leq \eta_k\leq\eta_i\} \;+\;
\#\{k<i|j\leq \eta_k+1\leq\eta_i\} \nonumber\\
{a'}(s)=j - 1 &&   \label{guion}
\end{eqnarray}
The following associated quantities occur frequently 
in the theory of the Macdonald polynomials. 
\begin{eqnarray}
d_{\eta}(q,t) &:=& \prod_{s\in\eta} 
\left( 1-q^{a(s)+1}t^{l(s)+1} \right)
\quad
d'_{\eta}(q,t):=\prod_{s\in\eta} 
\left( 1-q^{a(s)+1}t^{l(s)} \right)
\quad \nonumber\\
e_{\eta}(q,t) &:=& \prod_{s\in\eta} 
\left( 1-q^{a'(s)+1}t^{n-l'(s)} \right)
\quad
e'_{\eta}(q,t):=\prod_{s\in\eta} 
\left( 1-q^{a'(s)+1}t^{n-1-l'(s)} \right)
\quad \nonumber\\
b_{\eta}(q,t) &:=& \prod_{s\in\eta} 
\left( 1-q^{a'(s)}t^{n-l'(s)} \right)
\quad
l(\eta):= \sum_{s\in\eta}l(s) 
\quad \nonumber\\
l'(\eta)&:=& \sum_{s\in\eta}l'(s) \nonumber\\
\label{constants}
\end{eqnarray}
All these quantities are equal to one if $\eta=0$.
For future reference some properties of these 
quantities, easily derivable from \cite{Sa96}, are listed.
\begin{lemma} \label{lem2.1}
Let  $\Psi\eta:=(\eta_{2}\ldots\eta_{n},\eta_{1}+1)$
 and $\delta_{i, \eta}:={\bar{ \eta }}_{i}-{\bar{\eta}}_{i +1}$. 
We have
\begin{eqnarray}  \nonumber
\frac{d_{\Phi\eta }(q,t) }{ d_{\eta }(q,t) } =
\frac{e_{\Phi\eta }(q,t) }{ e_{\eta }(q,t) }=
1-qt^{n+{\bar{\eta}}_{1}}
&&
\frac{{d'}_{\Phi\eta }(q,t) }{ {d'}_{\eta }(q,t) } =
\frac{{e'}_{\Phi\eta }(q,t) }{ e'_{\eta }(q,t) }
=1-qt^{n-1+{\bar{\eta}}_{1}}
\\  \nonumber
l(\Phi\eta)=l(\eta)+\#\{k|k>1,\eta_{k} \leq \eta_1 \}  
&&
{l'}(\Phi\eta)={l'}(\eta)+n-1-\#\{k|k>1,\eta_{k} \leq \eta_1 \}
\end{eqnarray}
If $\eta_i > \eta_{i+1}$ we have
\begin{eqnarray} \nonumber
\frac{d_{s_{i} \eta }(q,t) }{ d_{\eta }(q,t) }
 =\frac{1-t^{\delta _{i, \eta}+1} }{ 1-t^{\delta _{i, \eta}} }
&&  \nonumber
\frac{{d'}_{s_{i} \eta }(q,t) }{ {d'}_{\eta }(q,t) }
 =\frac{1-t^{\delta _{i, \eta}} }{ 1-t^{\delta _{i, \eta}-1} }
\\ \nonumber
e_{s_{i}\eta} = e_{\eta}   && {e'}_{s_{i}\eta} = {e'}_{\eta}
\\  \nonumber
l(s_{i}\eta) = l(\eta)+1   && {l'}(s_{i}\eta) = {l'}(\eta)
\end{eqnarray}
\end{lemma}
The $q$-gamma function is defined by
\begin{equation}
\Gamma_{q}(x):=(q;q)_{x}(1-q)^{1-x}  \qquad
0<q<1
\end{equation}
where 
\begin{equation}
(q;q)_{a} :=\frac{ (q;q)_{\infty} }{ (q^{a};q)_{\infty} }
\qquad (b;q)_{\infty}:=\prod_{i=0}^{\infty}(1-bq^{i})
\end{equation}
We remark that with the generalised factorial defined by
\begin{eqnarray}
[q^{x}]_\eta^{(q,t)}&:=&\prod_{s\in\eta^{+}} 
\left( t^{{l'}(s)}-q^{a(s)+x} \right)
\nonumber\\
&=& t^{{l'}(\eta^{+})}(1-q)^{|\eta |}
\prod_{i=1}^n
\frac{\Gamma_{q}(x-\lambda (i-1) + \eta_{i}^{+}) 
}{ \Gamma_{q}(x-\lambda (i-1))} 
\label{cajamarca}
\end{eqnarray}
we have
\begin{equation}  \label{e-G}
e_{\eta}(q,t)=t^{-l'(\eta^+)}[q^{1+\lambda n}]_{\eta^{+}}^{q,t},\quad
{e'}_{\eta}(q,t)=t^{-l'(\eta^+)}
[q^{1+\lambda (n-1)}]_{\eta^{+}}^{q,t},\quad
b_{\eta}(q,t)=t^{-l'(\eta^+)}[q^{\lambda n}]_{\eta^{+}}^{q,t}
\end{equation}
Given a Laurent polynomial $f$ let CT$(f)$ denote the 
constant term in the laurent expansion of $f$ with 
respect to the variables $x_1,\ldots,x_{n}$. 
In the more general case where $f$ is not a 
laurent polynomial let CT$(f)$ denote the 
corresponding Fourier integral. The non-symmetric 
Macdonald polynomials have the following 
orthogonality property, which can be deduced from (\ref{E2}).
\begin{prop} {\rm\cite{Ma95}}
Given any two polynomials $f(x;q,t) $ and $g(x;q,t) $ define
the scalar product
\begin{equation}
\langle f,g \rangle_{q,t} :={\rm CT}\left(f(x;q,t) g(\frac{1}{ x}; 
\frac{1}{ q}, \frac{1}{ t}) W(x) \right)
\end{equation}
where
\begin{equation}
W(x):=W(x;q,t):=\prod_{1\leq i<j \leq n} 
(\frac{x_{i} }{ x_{j}} ;q)_{\lambda}    
(q\frac{x_{j} }{ x_{i}} ;q)_{\lambda}  
\end{equation}
The polynomials $E_{\eta}(x;q,t) $  
form an orthogonal set with respect to  
$\langle \cdot ,\cdot \rangle_{q,t}$.
\end {prop} 

A consequence of this is that the non-symmetric
Macdonald polynomials are able to be constructed
by means of a Gram-Schmidt procedure. Let 
$\eta^{(1)}\prec \cdots \prec \eta^{(p)}$
be a chain of compositions satisfying
\begin{equation} \label{complete-chain}
\mbox{If}\quad \eta^{(i)}\prec\mu\prec\eta^{(i+1)}\quad
 \mbox{then $\mu = \eta^{(i)}$ or $\mu = \eta^{(i+1)}$}
\end{equation}
The non-symmetric Macdonald polynomial $E_{\eta^{(p)}}$
can be determined as the unique polynomial satisfying
(\ref{E1}) which is orthogonal to all $E_{\eta^{(i)}}$
with $i<p$.

The non-symmetric Macdonald polynomials are
elements of the ring of $n$ variable polynomials
whose coefficients are elements of the field
$\mathbb{Q}(q,t)$ of rational functions in $q$ and $t$.
As in the symmetric case the symbols
$q$ and $t$ can be interpretated as indeterminants
or as parameters ranging over
$0<q,\, t<1$. 
Let the hat symbol \, $\hat{\,}$ \, denote the involution on this ring 
 which sends $x_{i}\mapsto x_{n-i+1}$,
$q \mapsto q^{-1}$ and $t \mapsto t^{-1}$.
Extend this operator to act on operators so that
for any operator $T$ and polynomial $f$, $\hat{T} \hat{f}= (\widehat{Tf})$.
We define the dual non-symmetric Macdonald polynomial by
${\hat{E}}_{\eta}(x;q,t):=E_{\eta}(\underline{x};q^{-1},t^{-1})$
where $\underline{x}:=(x_{n}, \ldots,x_{1})$.

These polynomials are uniquely determined by the conditions
\begin{eqnarray}  \label{dual}
{\hat{E}}_{\eta}(x;q,t) &=&
x^{\underline{\eta}} + \sum_{\nu \acute{\prec} \underline{\eta}} 
c_{\eta\nu} x^{\nu} \\ \label{dual.b}
{\hat{Y}}_{i} \, {\hat{E}}_{\eta}(x;q,t) &=& 
t^{-\bar{\eta}_i}\,{\hat{E}}_{\eta}(x;q,t) \qquad
1\leq i \leq n
\end{eqnarray}
where $\underline{\eta}:=(\eta_n ,\ldots \eta_1)$ and
the partial order $ \acute{\prec}$ is defined on compositions
so that
\begin{equation}
\nu\acute{\prec}\eta \quad\mbox{iff}\quad \nu^+ < \eta^+
\quad\mbox{or in the case $\nu^+ = \eta^+$}\quad \nu > \eta
\end{equation}
Note that if $\nu^{+}=\eta^{+}$ then
$ \nu\acute{\prec}\eta$ iff  $ \nu\succ\eta$.

The dual non-symmetric Macdonald polynomials are simply related
to the non-symmetric  Macdonald polynomials by means of the
Demazure-Lustig operators.
\begin{lemma}  \label{E-T-hatE}
\begin{eqnarray}   \label{E-T-hatE.1}
\mbox{a)}\qquad T_{(n,\ldots,1)}\,{\hat{E}}_{\eta}(x;q,t)&=&
t^{\#\{(i,j)|i<j,\eta_{i}\geq\eta_{j}\}}E_{\eta}(x;q,t)
\\
\mbox{b)}\qquad T_{(n,\ldots,1)}\,E_{\eta}(x;q,t)&=&
t^{\#\{(i,j)|i<j,\eta_{i}\leq\eta_{j}\}}{\hat{E}}_{\eta}(x;q,t)
\end{eqnarray}
\end{lemma}

\noindent {\it Proof.} \quad 
We shall only consider (a) as the proof of (b) is similar.
A direct calculation
using (\ref{t.action}) reveals that
\begin{equation} 
T_{(n,\ldots,1)}\,x^{\eta}=t^{\#\{(i,j)|i<j,\eta_{i}\leq \eta_{j}\}}
x^{\underline{\eta}}+\sum_{\mu\prec \underline{\eta}}a_{\mu} x^{\mu}
\end{equation}
It follows that 
\begin{equation}  \label{E-T-hatE.2}
T_{(n,\ldots,1)}\,{\hat{E}}_{\eta}(x;q,t)=
t^{\#\{(i,j)|i<j,\eta_{i}\geq\eta_{j}\}}
\left( x^{\eta}+
\sum_{\mu\prec \eta}a'_{\mu} x^{\mu} \right)
\end{equation}
It suffices then to show that given a chain 
$\eta^{(1)}\prec \cdots \prec \eta^{(p)}=\eta$
satisfying (\ref{complete-chain}),
$T_{(n,\ldots,1)}\,{\hat{E}}_{\eta^{(p)}}$ is orthogonal
to $E_{\eta^{(i)}}$ for all $i<p$.
This will be done by induction.
If $\mu$ is a minimal composition 
under the partial ordering $\prec$
then $E_{\mu}(x;q,t)=x^{\mu}$.
It follows from (\ref{E-T-hatE.2}) that (\ref{E-T-hatE.1}) 
is true for the composition $\eta^{(1)}$.
Suppose (\ref{E-T-hatE.1}) is true for 
$\eta^{(1)},\ldots ,\eta^{(r-1)}$. Then for any
$k<r$
\begin{equation}
\langle E_{\eta^{(k)}},T_{(n,\ldots,1)}\, {\hat{E}}_{\eta^{(r)}}
\rangle_{q,t}
=t^{-\#\{(i,j)|i<j,\eta_{i}\geq\eta_{j}\}}
\langle T_{(n,\ldots,1)}\, {\hat{E}}_{\eta^{(k)}},T_{(n,\ldots,1)}\,
{\hat{E}}_{\eta^{(r)}}
\rangle_{q,t}
\end{equation}
Since $T^{-1}_{i}$ is the adjoint operator of $T_{i}$ 
 with respect to
$\langle\cdot ,\cdot \rangle_{q,t}$ \cite{Ma95a}
and $\langle\hat{f} ,\hat{g} \rangle_{q,t}=
\langle f ,g \rangle_{q,t}$ we have
\begin{eqnarray} \nonumber
\langle E_{\eta^{(k)}},T_{(n,\ldots,1)}\, {\hat{E}}_{\eta^{(r)}}
\rangle_{q,t}
&=&t^{-\#\{(i,j)|i<j,\eta_{i}\geq\eta_{j}\}}
\langle  {\hat{E}}_{\eta^{(k)}},
{\hat{E}}_{\eta^{(r)}} \rangle_{q,t}
\\  \nonumber 
&=&
t^{-\#\{(i,j)|i<j,\eta_{i}\geq\eta_{j}\}}
\langle  E_{\eta^{(k)}},
E_{\eta^{(r)}} \rangle_{q,t}
\\
&=& 0
\end{eqnarray}
\hfill $\Box$

Next we revise the construction of the 
symmetric and $q$-antisymmetric Macdonald 
polynomials from the non-symmetric Macdonald 
polynomials.
This requires introducing $q$-analogues of the symmetrization 
and antisymmetrization operators defined by \cite{Ma95a}
\begin{equation}
U^{+}:=\sum_{\sigma \in S_{n}}T_{\sigma}  \qquad
U^{-}:=\sum_{\sigma \in S_{n}}(-t)^{-l(\sigma)}T_{\sigma}    
\end{equation}
where  $l(\sigma):=\#\{(i,j)| i<j, {\sigma}_{i} >{\sigma}_{j}\} $ 
is the length of the permutation $\sigma$.
These operators have the following properties 
\begin{eqnarray}  \label{TU+}
&& T_{i}^{\pm}U^{+} = U^{+} T_{i}^{\pm}=t^{\pm} U^{+}
\\ \label{TU-}
&& T_{i}^{\pm}U^{-} = U^{-} T_{i}^{\pm}=- U^{-}
\end{eqnarray}
From these properties it can be deduced that
\begin{equation}  \label{Us=-U}
U^{-}s_{i}=-U^{-}
\end{equation}
from which it follows that
\begin{equation}  \label{U=0}
U^{-} x^{\mu} =0 \qquad\mbox{if $\mu_{i}=\mu_{j}$ for
$i\neq j$}
\end{equation}
For $\delta := (n-1,\ldots,1,0)$ it is also the case that
\begin{equation} \label{Ux}
t^{n(n-1)/2} U^{-} x^{\delta}=\prod_{1\leq i<j \leq n}
(tx_{i}-x_{j}):=\Delta_{t}(x)
\end{equation}
Now, when acting on symmetric functions, the Macdonald operator
$ D_{n}^{1}(q,t) $ can be decomposed in terms of the 
Cherednik operators  according to \cite{KN}
\begin{equation}
D_{n}^{1}(q,t)=t^{n-1}\sum_{i=1}^{n}Y_{i}
\end{equation}
Since the operator $U^{+}$ commutes with $\sum_{i=1}^{n}Y_{i}$
it follows from (\ref{E1}) and (\ref{E2}) that there exist unique 
symmetric polynomials indexed by partitions which satisfy
\begin{eqnarray} \label{Pdef1}
\left(\sum_{i=1}^{n}Y_{i} \right) P_{\ka}(x;q,t)  
&=&\left(\sum_{i=1}^{n}t^{{\bar{\ka}}_{i}} \right) P_{\ka}(x;q,t)  
\\ \label{Pdef2}
P_{\ka}(x;q,t) &=&m_{\ka}(x) +\sum_{\mu<\ka}u_{\ka\mu} m_{\mu}(x)
\end{eqnarray}
From Section 1 these are the symmetric Macdonald polynomials.
One has the relation
\begin{equation} \label{U+E}
P_{\eta^{+}}(x;q,t) =\frac{1}{ \alpha_{\eta}(q,t)}U^{+}E_{\eta}(x;q,t)
\end{equation}
for scalars $\alpha_{\eta}(q,t)$.
We can also define the $q$-antisymmetric Macdonald polynomials 
\cite{Ma95a}. 
The $q$-antisymmetric monomial ${m'}_{\ka} $, indexed by the 
partition $\ka$ with non-repeating parts, is
\begin{equation}
{m'}_{\ka}:=U^{-}x^{\ka}
\end{equation}
A function $f$ is $q$-antisymmetric if for all $i$, $T_{i}f=-f$.
The $q$-antisymmetric monomials are a basis for the 
analytic $q$-antisymmetric functions. The $q$-antisymmetric 
Macdonald polynomials $S_{\ka}(x;q,t) $ are indexed by partitions 
with non-repeating parts and can be defined by the following conditions.
\begin{eqnarray} \label{S1}
\left(\sum_{i=1}^{n}Y_{i} \right) S_{\ka}(x;q,t)  
&=&\left(\sum_{i=1}^{n}t^{{\bar{\ka}}_{i}} \right) S_{\ka}(x;q,t)  
\\ \label{S2}
S_{\ka}(x;q,t) &=&{m'}_{\ka}(x) +\sum_{\mu<\ka}v_{\ka\mu} {m'}_{\mu}(x)
\end{eqnarray}
Analogous to the derivation of (\ref{U+E}) we have
\begin{equation}  \label{Us}
S_{\eta^{+}}(x;q,t) =\frac{1}{ \beta_{\eta}(q,t)}U^{-}E_{\eta}(x;q,t)
\end{equation}
The symmetric and $q$-antisymmetric Macdonald polynomials 
can also be expressed as linear combinations of the non-symmetric 
Macdonald polynomials.
\begin{lemma} {\rm\cite{Ma95a}}\label{lemP-E}
\begin{eqnarray}
\mbox{a)}\qquad \label{P-E}
P_{\ka}(x;q,t) &=&\sum_{\eta :\eta^{+}=\ka}
\frac{{d'}_{\eta ^{+}}(q,t) }{ {d'}_{\eta }(q,t)  }
E_{\eta}(x;q,t)
\\
\mbox{b)}\qquad \label{S-E}
S_{\ka}(x;q,t) &=&\sum_{\sigma\in S_{n} }{(-t)}^{-l(\sigma)}
\frac{{d}_{\sigma (\ka)}(q,t) }{ {d}_{\ka }(q,t)  }
E_{ \sigma (\ka)}(x;q,t)
\end{eqnarray}
\end{lemma}
{\it Proof.} \quad A simple generalisation of the derivation 
of the analogue results in the case of the Jack polynomials \cite{BDF}.
\hfill $\Box$ \medskip

It immediately follows from the orthogonality of the non-symmetric 
Macdonald polynomials and Lemma \ref{lemP-E} that
\begin{prop} {\rm\cite{Ma95a}}
 Both the symmetric Macdonald polynomials
$\{ P_{\ka}(x;q,t) \}$ and the $q$-antisymmetric Macdonald 
polynomials $\{ S_{\ka}(x;q,t) \}$ form orthogonal sets with 
respect to  $\langle \cdot ,\cdot \rangle_{q,t}$.
\end{prop}

It follows that both the symmetric and
$q$-antisymmetric Macdonald 
polynomials are able to be constructed by means of a
Gram-Schmidt procedure similar to that in the case of the 
non-symmetric polynomials.

The dual non-symmetric Macdonald polynomials share many
properties with the non-symmetric Macdonald polynomials.
In particular they are equally able to serve as building blocks
for the symmetric and $q$-antisymmetric Macdonald polynomials.
This is explained by the following results.
\begin{lemma} \label{U-hatU}
\begin{eqnarray}
\mbox{a)} &\qquad &{\hat{U}}^{+}=t^{-\frac{n(n-1) }{ 2}}U^{+}
\\
\mbox{b)} &\qquad &{\hat{U}}^{-}=t^{\frac{n(n-1) }{ 2}}U^{-}
\end{eqnarray}
\end{lemma}
\noindent {\it Proof.} \quad 
It is well known \cite{We88} that given a Hecke algebra
${\cal H}_{n}(t)$ generated by $T_1$, $\ldots$,$T_{n-1}$,
with $t$ not a root of unity,
there exist unique elements (up to scalar multiplication)
$\alpha,\beta \in {\cal H}_{n}(t) $  such that for all $i=1, \ldots,n-1$,
\begin{equation}
T_{i}\,\alpha=t\alpha \quad \mbox{and} \quad T_{i}\,\beta=-\beta
\end{equation}
It follows from the definitions that ${\hat{T}}_{i}=T^{-1}_{n-i}$.
 ${\hat{U}}^{+}$ and ${\hat{U}}^{-}$ are  then both
elements of ${\cal H}_{n}(t)$.
Using (\ref{TU+}) and (\ref{TU-}) we have 
\begin{eqnarray}
&&T_{i}{\hat{U}}^{+}=(\widehat{\hat{T_{i}}U^{+}})=
(\widehat{{T}^{-1}_{n-i}U^{+}})
=t{\hat{U}}^{+}
\\ &&
T_{i}{\hat{U}}^{-}=(\widehat{\hat{T_{i}}U^{-}})
=(\widehat{{T}^{-1}_{n-i}U^{-}})
=-{\hat{U}}^{-}
\end{eqnarray}
Hence ${\hat{U}}^{+}=c_{1}U^{+}$ and ${\hat{U}}^{-}=c_{2}U^{-}$.
 Equating
coefficients of the basis $\{ T_{\sigma}: \sigma \in S_{n} \}$
of ${\cal H}_{n}(t)$
reveals the coefficients $c_{1}$ and $c_{2}$ to be $t^{-n(n-1)/2}$
and $t^{n(n-1)/2}$ respectively.
\hfill $\Box$ \medskip

\begin{lemma} \label{L4.1}
\begin{eqnarray}
\mbox{a)} &\qquad& P_{\kappa}(x;q,t)=
{\hat{P}}_{\ka}(x;q,t)=P_{\kappa}(x;q^{-1},t^{-1}) \label{4.1.1}
\\
\mbox{b)}  &\qquad& S_{\ka}(x;q,t)=
(-t)^{-\frac{n(n-1)}{ 2}}
{\hat{S}}_{\ka}(x;q,t)
\end{eqnarray}
\end{lemma}
\noindent {\it Proof.} \quad 
We shall consider only the second identity as
(a) is well known and is proven in a similar way as (b).
It follows from Lemma \ref{U-hatU}
that
${\hat{m'}}_{\lambda}(x)=(-t)^{n(n-1)/2}{m'}_{\lambda}(x)$.
Using the defining property (\ref{S2}) we then have
\begin{equation}
{\hat{S}}_{\ka}(x;q,t)=(-t)^{\frac{n(n-1)}{ 2}}
\left( {m'}_{\ka}(x)+ \sum_{\mu < \ka}{\hat{v}}_{\ka \mu}
{m'}_{\mu}(x) \right)
\end{equation}
Since $\{ (-t)^{-n(n-1)/ 2} {\hat{S}}_{\ka}(x;q,t) \}$ 
is orthogonal with respect to
$\langle \cdot,\cdot \rangle_{q,t}$ and posseses the 
triangular structure (\ref{S2}) (b) must be true.
\hfill $\Box$ \medskip

Using the above two lemmas in conjunction with (\ref{U+E}),
 (\ref{Us}) and lemma \ref{lemP-E}
we obtain the following two lemmas.
\begin{lemma}
\begin{equation}
\mbox{a)} \qquad P_{\eta^+}(x;q,t)=
\frac{t^{-\frac{n(n-1)}{ 2}}}{ \alpha_{\eta}(q^{-1},t^{-1})}
U^{+}{\hat{E}}_{\eta}(x;q,t)
\end{equation}
\begin{equation}
\mbox{b)} \qquad S_{\eta+}(x;q,t)=
\frac{(-1)^{n(n-1)/2} }{ \beta_{\eta}(q^{-1},t^{-1})}
U^{-}{\hat{E}}_{\eta}(x;q,t)
\end{equation}
\end{lemma}
\begin{lemma}
\begin{eqnarray}  \label{dsymsum}
\mbox{a)} & & P_{\kappa}(x;q,t)=
\sum_{\eta: \eta+=\ka} \frac{{d'}_{\eta+}(q^{-1},t^{-1}) }{
{d'}_{\eta}(q^{-1},t^{-1})}
{\hat{E}}_{\eta}(x;q,t)
\\
\mbox{b)}  & & S_{\ka}(x;q,t)=
\sum_{\sigma \in S_{n}}(-t)^{l(\underline{\sigma})}
\frac{d_{\sigma(\ka)}(q^{-1},t^{-1}) 
}{
d_{\ka}(q^{-1},t^{-1})}
{\hat{E}}_{\sigma(\ka)}(x;q,t)
\end{eqnarray}
where for any permutation $\sigma$,
$\underline{\sigma} :=(\sigma_n ,\ldots ,\sigma_1)$.
\end{lemma}

\setcounter{equation}{0}
\section{Non-symmetric Macdonald Polynomial Theory}
In this section we will derive some of the basic properties
of the non-symmetric Macdonald polynomials independently of
the theory of the symmetric Macdonald polynomials.
 A required preliminary result is the
Cauchy type formula for the non-symmetric Macdonald polynomials.
\begin{prop} {\rm\cite{MN}}\label{E-sum}
\begin{equation}  \label{E-sum1}
\Omega (x,y;q,t)= \sum_{\eta}\frac{1}{ u_{\eta}(q,t) } E_{\eta}(x;q,t)
E_{\eta}(y;q^{-1} ,t^{-1}), \quad
\Omega (x,y;q,t):= \prod_{i=1}^{n}\frac{1}{ (x_{i}y_i ;q)_{\lambda +1} }
\prod_{1\leq i<j \leq n}\frac{1}{ (x_{i}y_j ;q)_{\lambda } 
(qx_{j}y_i ;q)_{\lambda } }
\end{equation}
\end{prop}
\noindent  {\it Remark.} \quad
Define the lengh of a composition to be 
${\rm length}(\eta):={\rm max}\{i|\eta_{i}\neq 0\}$.
The non-symmetric Macdonald polynomials have the following stability
property.
\begin{equation}  \label{E-stability}
E_{\eta}(x_1 ,\ldots, x_{n-1},0;q,t) =\left\{ \begin{array}{ll}
E_{\eta}(x_1 ,\ldots, x_{n-1};q,t)\quad & {\rm length}(\eta)\leq n-1
\\
0 & {\rm else}
\end{array}
\right. \end{equation}
Applying this property to (\ref{E-sum1}) shows that 
the scalars $ u_{\eta}(q,t) $ are 
independent of $n$.\medskip

Dunkl has introduced a family of multivariable polynomials which 
allow a workable treatment of some important constructions and 
has a close relationship to the theory of the non-symmetric 
Macdonald polynomials. The $q$-analogue of these polynomials 
are the polynomials $q_{\eta}(x;q,t) $ defined by
\begin{equation}
\Omega (x,y;q,t) :=\sum_{\eta} q_{\eta}(x;q,t) y^{\eta}
\end{equation}
\begin{cor}  \label{<>q1}
Define a scalar product by $\langle E_{\nu}(x;q,t)  , 
\langle E_{\nu}(x;q^{-1},t^{-1})  \rangle_{q } 
:= u_{\eta}(q,t) \delta_{\nu\eta} $. We have
\begin{equation}
\langle q_{\nu }(x;q,t), x^{\eta}\rangle_{q }
=\delta_{\nu\eta}
\end{equation}
Hence the $ q_{\eta}(x;q,t) $ are a basis for the 
multivariable polynomials with coefficients in $\mathbb{Q}(q,t)$.
\end{cor}
{\it Proof.} \quad From the triangular structure of the 
non-symmetric Macdonald polynomials, 
$\{ \frac{1}{ u_{\eta}(q,t) } E_{\eta}(x;q,t)\}$ and 
$ \{ E_{\eta}(x;q^{-1} ,t^{-1}) \}$ are basis
for the multivariable polynomials. The scalar product 
$\langle \cdot  , \cdot  \rangle_{q } $ is then well defined.
An argument similar to MacdonaldÕs 
\cite[p310-11]{Ma95} can now be used to 
show that (\ref{<>q1}) is equivalent to (\ref{E-sum1}).
\hfill $\Box$ \medskip

The non-symmetric Macdonald polynomials can be computed recursively be just 
two kinds of operators. The first are the Demazure-Lustig operators
 $T_{i} $, $1\leq i\leq n-1$. The
second, introduced by Baker and Forrester \cite{BF97g}, 
is the raising-type
operator 
\begin{equation}
\Phi_{q} :=x_{n} T_{n-1}^{-1} \cdots T_{2}^{-1}\,T_{1 }^{-1} 
=t^{i-n} T_{n-1} \cdots T_{i}\, x_{i}\,
T_{i-1}^{-1} \cdots T_{1 }^{-1} 
\end{equation}
These operators have the following action on the non-symmetric
Macdonald polynomials \cite{BF97g},\cite{MN}
\begin{equation}  \label{Phi-E}
\Phi_{q} E_{\eta}(x;q,t) =t^{-\#\{ i|i>1, \eta_{i} \leq \eta_{1} \} }
E_{\Phi\eta}(x;q,t)
\end{equation}
and
\begin{equation} \label{TE}
T_i\,E_{\eta} = \left\{ \begin{array}{ll}
\left( \frac{t-1}{ 1-t^{-\de_{i\eta}}}\right)\,E_{\eta}
+ t\:E_{s_i\eta} & \eta_i < \eta_{i+1} \\
t\:E_{\eta} & \eta_i = \eta_{i+1} \\
\left (\frac{t-1}{ 1-t^{-\de_{i\eta}}}\right)\,E_{\eta} +
\frac{(1-t^{\de_{i\eta}+1})(1-t^{\de_{i\eta}-1})}{
(1-t^{\de_{i\eta}})^2}\,E_{s_i\eta} & \eta_i > \eta_{i+1}
\end{array}\right. 
\end{equation}
Using these operators, it is simple to derive the following two identities by verifying that the respective quantities satisfy the same recursion relationships.
\begin{prop} {\rm\cite{Ch95a}}\label{propE(t)}
 Let $ t^{\underline{\delta}} :=(1,t,\ldots,t^{n-1})$.
We have
\begin{equation} \label{E(t)}
E_{\eta}( t^{\underline{\delta}} ;q,t) =t^{l(\eta)}
\frac{e_{\eta}(q,t) }{ d_{\eta}(q,t) }
\end{equation}
\end{prop}

\noindent {\it Proof. }  \quad Noting that for any function $f=f(x)$
\begin{equation}  \label{Tf}
(T_{i}f)( t^{\underline{\delta}})=tf(t^{\underline{\delta}})
\end{equation}
shows that 
\begin{eqnarray}  \nonumber
\left. \left( \Phi_{q} E_{\eta}( x ;q,t) \right) 
\right|_{x= t^{\underline{\delta}} }
&=&
\left. \left(t^{1-n} T_{n-1} \cdots T_{1}x_{1}E_{\eta}( x ;q,t) \right) 
\right|_{x= t^{\underline{\delta}} }
\\
&=& E_{\eta}( t^{\underline{\delta}} ;q,t)
\end{eqnarray}
Using (\ref{Phi-E}) we then obtain
\begin{equation}  \label{3.3.1}
E_{\Phi\eta} ( t^{\underline{\delta}} ;q,t) =
t^{ \#\{ i|i>1, \eta_{i} \leq \eta_{1} \} }
E_{\eta}( t^{\underline{\delta}} ;q,t)
\end{equation}
Supposing $\eta_{i }>\eta_{i +1}$ and applying (\ref{Tf}) to 
(\ref{TE}) and rearranging, we also obtain
\begin{equation} \label{3.3.2}
E_{s_i \eta} ( t^{\underline{\delta}} ;q,t) =
t\frac{1-t^{\delta_{i,\eta}} }{ 1-t^{\delta_{i,\eta}+1}     }
E_{\eta}( t^{\underline{\delta}} ;q,t)
\end{equation}
The relations (\ref{3.3.1}) and (\ref{3.3.2}) uniquely determine
$ E_{\eta}( t^{\underline{\delta}} ;q,t)$ given $ E_{0}
( t^{\underline{\delta}} ;q,t)$. Since Proposition \ref{propE(t)} is 
obviously true for the case $\eta=0$ all that remains is to show that
 the right hand side of (\ref{E(t)}), RHS$(\eta)$ say, obeys these relations. Using Lemma
\ref{lem2.1} we have
\begin{equation}
 \frac{{\rm RHS} (\Phi\eta) }{{\rm RHS} (\eta) }=
t^{ \#\{ i|i>1, \eta_{i} \leq \eta_{1} \} }
\end{equation}
While supposing $\eta_{i}>\eta_{i+1}$ and again using Lemma
\ref{lem2.1} we have
\begin{equation}
{\rm RHS} (s_{i}\eta) =t
\frac{1-t^{\delta_{i,\eta}} }{ 1-t^{\delta_{i,\eta}+1}     }
{\rm RHS} (\eta)
\end{equation}
\hfill $\Box$ \medskip

\begin{prop}   \label{NE}
Write ${\cal N}_{\eta}^{(E)}(q,t) :=\langle E_{\eta} , E_{\eta}\rangle_{q,t}$. We have
\begin{equation}  \label{NE.1}
\frac{{\cal N}_{\eta}^{(E)} (q,t) }{ {\cal N}_{0}^{(E)}(q,t) }=
\frac{{d'}_{\eta} (q,t) {e}_{\eta} (q,t) }{
{d}_{\eta} (q,t) {e'}_{\eta}(q,t)}
\end{equation}
\end{prop}
\noindent {\it Remark.}
Macdonald \cite{Ma95a} and Cherednik \cite{Ch95a}
have derived (\ref{NE.1}) although in a different form.
\medskip 

\noindent {\it Proof.} \quad Using (\ref{Phi-E}) we have
\begin{eqnarray}  \nonumber
\langle E_{\Phi\eta} , E_{\Phi\eta}\rangle_{q,t}
&=&\langle t^{\#\{ i|i>1, \eta_{i} \leq \eta_{1} \} } \Phi_{q}E_{ \eta} , 
t^{\#\{ i|i>1, \eta_{i} \leq \eta_{1} \} }\Phi_{q}E_{ \eta} \rangle_{q,t}
\\  \nonumber
&=&{\rm CT}\left( \left. t^{\#\{ i|i>1, \eta_{i} \leq \eta_{1} \} } x_n
\left( T_{n-1}^{-1} \cdots T_{1}^{-1} E_{ \eta} (x;q,t) \right)
\right. \right. 
\\   \nonumber
&& \qquad \qquad \times t^{-\#\{ i|i>1, \eta_{i} \leq \eta_{1} \} } x_{n}^{-1} \left. \left.
\left(  T_{n-1}^{-1} \cdots T_{1}^{-1} E_{ \eta}(x;q,t) 
\right) \right|_{\frac{1 }{ x} , \frac{1 }{ q} , \frac{1 }{ t}  }
W(x;q,t)\right)
\\  \nonumber
&=&
\langle T_{n-1}^{-1} \cdots T_{1}^{-1} E_{ \eta}(x;q,t), 
T_{n-1}^{-1} \cdots T_{1}^{-1} E_{ \eta}(x;q,t) \rangle_{q,t}
\\
&=&
\langle E_{\eta} , E_{\eta}\rangle_{q,t}
\end{eqnarray}
In the last line we have used the fact that $T_{i}^{-1}$ is the 
adjoint operator of $ T_{i} $ with respect to 
$\langle \cdot , \cdot \rangle_{q,t} $ \cite{Ma95a}.

Supposing $\eta_{i}< \eta_{i+1}$ and using (\ref{TE}) we have 
\begin{eqnarray}   \nonumber
\langle E_{ s_{i} \eta} , E_{ s_{i} \eta}\rangle_{q,t}
&=&\langle t^{-1}T_{i}E_{\eta} -\frac{1-t^{-1} }{ 1-t^{-\delta_{i,\eta}} }
E_{\eta} , 
t^{-1}T_{i}E_{\eta} -\frac{1-t^{-1} }{ 1-t^{-\delta_{i,\eta} } }
E_{\eta}\rangle_{q,t}
\\ \nonumber
&=&  \langle T_{i}E_{ \eta} , T_{i}E_{ \eta} \rangle_{q,t}
-t^{-1}\frac{1-t }{ 1-t^{\delta_{i,\eta}} }
\langle T_{i}E_{ \eta} , E_{ \eta} \rangle_{q,t}
\\
&& \qquad\qquad -t\frac{1-t^{-1} }{ 1-t^{-\delta_{i,\eta}}}
\langle E_{ \eta} , T_{i}E_{ \eta} \rangle_{q,t}
+ \frac{(1-t) (1-t^{-1}) }{  (1-t^{\delta_{i,\eta}}) 
(1-t^{-\delta_{i,\eta}})     }
\langle E_{ \eta} , E_{ \eta} \rangle_{q,t}  \label{NE.2}
\end{eqnarray}
Consider the right hand side of this expression. The first term simplifies by
again using the fact that that $T_{i}^{-1}$ and $ T_{i} $ are
adjoint operators, while the second and third terms simplify
by making further use of (\ref{TE}) and then noting that
for $\eta_{i}\neq \eta_{i+1}$, $E_{\eta}$ and $E_{s_{i} \eta}$
are orthogonal. We obtain after rearranging
\begin{equation}
\langle E_{ s_{i} \eta} , E_{ s_{i} \eta}\rangle_{q,t} =
\frac{(1-t^{\delta_{i,\eta}+1}) (1-t^{\delta_{i,\eta}-1})  }{
(1-t^{\delta_{i,\eta}})^{2} }
\langle E_{ \eta} , E_{ \eta}\rangle_{q,t}
\end{equation}
By replacing $\eta$ by $ s_{i} \eta$ and 
noting that if $\eta_i \neq \eta_{i+1}$ 
$\delta_{i,s_{i}\eta} =-\delta_{i,\eta} $ we see 
that in the case $\eta_{i}> \eta_{i+1}$
\begin{equation}  \label{NE.3}
\langle E_{ s_{i} \eta} , E_{ s_{i} \eta}\rangle_{q,t} =
\frac{(1-t^{\delta_{i,\eta}})^{2} }{
(1-t^{\delta_{i,\eta}+1}) (1-t^{\delta_{i,\eta}-1})  }
\langle E_{ \eta} , E_{ \eta}\rangle_{q,t}
\end{equation}
Using Lemma \ref{lem2.1} it is clear that the right hand
side of (\ref{NE.1}) satisfies both the recursion relations
(\ref{NE.2}) and (\ref{NE.3}). Since (\ref{NE.1}) is true in the 
trivial case $\eta=0$ Proposition \ref{NE} is true by induction.
\hfill $\Box$ \medskip

We shall now show that the multivariable $q$-binomial
theorem involving the non-symmetric Macdonald polynomials
can be deduced using  Propositions \ref{E-sum}
and \ref{propE(t)}.

\begin{prop} {\rm \cite{MN} } \label{Bi-E}
\begin{equation}  \label{Bi-E.1}
\prod_{i=1}^{n}\frac{1}{ (x_{i}; q)_{r}}=
\sum_{\eta} \frac{[q^{r}]_{\eta^+} }{ u_{\eta}(q,t) d_{\eta}(q,t) }
E_{\eta}(x;q,t)
\end{equation}
\end{prop} 
\noindent {\it Remark.} \quad  The expression on the right hand
side of (\ref{Bi-E.1}) will be able to
 be simplified using (\ref{u-eta.1}).
\medskip

\noindent {\it Proof.} \quad
In (\ref{E-sum1}) first replace $n$ by $kn$ for some 
$k\in {\mathbb{Z}}_{>0} $
and then substitute $y_{j}=t^{kn-j}$ and let 
$x_{n+1}= \cdots =x_{kn}=0 $.  Since 
$E_{\eta}(cx) = c^{|\eta|}E_{\eta}(x) $
we can use Proposition \ref{propE(t)} to obtain
\begin{equation}
\prod_{i=1}^{n}\frac{1}{ (x_{i}; q)_{kn\lambda +1}}=
\sum_{\eta}
\frac{\left. t^{(kn-1)|\eta|-l(\eta)} e_{\eta}(\frac{1}{ q},\frac{1}{ t}) 
\right|_{n\rightarrow kn} }{
u_{\eta}(q,t) d_{\eta}(\frac{1}{ q},\frac{1}{ t}) }
E_{\eta}(x_{1} ,\ldots, x_{n},0,\ldots ,0;q,t)
\end{equation}
Making use of (\ref{e-G}) ,   Lemma \ref{lem2.1}, the stability
property (\ref{E-stability}) and the identity
\begin{equation} \label{e/d(1/q,1/t)}
\frac{ \left. e_{\eta}(q^{-1} ,t^{-1})\right|_{n\rightarrow kn}  
}{ d_{\eta}(q^{-1} ,t^{-1}) }
=t^{\left( l(\eta)+{l'}(\eta) -(kn-1)|\eta|\right) }
\frac{ \left. e_{\eta}(q ,t)\right|_{n\rightarrow kn}  
}{ d_{\eta}(q ,t) }
\end{equation}
we obtain for $k\in {\mathbb{Z}}_{>0}$
\begin{equation}  \label{Bi-E.2}
\prod_{i=1}^{n}\frac{1}{ (x_{i}; q)_{kn\lambda +1}} =
\sum_{\eta} \frac{[q^{ kn\lambda +1}]_{\eta^+} }{
u_{\eta}(q,t) d_{\eta}(q,t) }
E_{\eta}(x;q,t)
\end{equation}

To show that (\ref{Bi-E.2}) is true for all $k\in \mathbb{R}$
we require $ u_{\eta}(q,t) $ to be able to be written as
a power series in $q$ and $t$ for all $0<q,\, t<1$. 
This falls out of the proof of Proposition \ref{E-sum}
in \cite{MN} by using expansions \cite[(2.3),(2.4)]{BF} and by noting
that the coefficients of the non-symmetric Macdonald
polynomials can be written as 
power series in $q$ and $t$ for all $0<q,\, t<1$.
Both sides of (\ref{Bi-E.2}) are then power series in
$x,q,t$ and $t^{k}$. Equating the coefficients with respect
to $q$ and $x$ we can apply the following lemma to show
that the $q$-binomial theorem (\ref{Bi-E.1}) is true
for all $k\in \mathbb{R}$.
\hfill $\Box$ \medskip

\begin{lemma}{\rm\cite{St88}} {\rm\label{Stem}}
Let $F(z,q)$ and $G(z,q)$ be formal power
series in $z$ and $q$. If $F(q^{k} ,q)=G(q^{k},q)$ for infinitely
many integers $k\geq 0$ then $F=G$.
\end{lemma}

\setcounter{equation}{0}
\section{\mbox{A Generalisation of the $q$-Selberg Integral}}
The $q$-Selberg Integral, as formulated by Askey \cite{As80} and 
subsequently
proved by Kadell (\cite{Ka88} and Habseiger \cite{Ha},
 has been extended by Kadell
\cite{Ka88} and Kaneko \cite{Kan96} to involve the symmetric Macdonald polynomial
as a factor in the integrand. An equivalent formulation
of this result is as a constant term identity which generalises
the $q$-Morris identity \cite{Kan97}. Here this result
will itself be extended in that the symmetric
Macdonald polynomial will be replaced by the 
non-symmetric Macdonald polynomial. The derivation of
this identity will also yield a new derivation of the
$q$-Selberg integral as well as allowing us to specify the constant
$u_{\eta}(q,t)$ appearing in (\ref{E-sum1}).
The derivation 
is based on the multivariable
$q$-binomial theorem (\ref{Bi-E.1}). 

Since $\{ E_{\eta}(x;q,t) \}$
is an orthogonal basis for multivariable analytic functions  with respect
to $\langle \cdot,\cdot \rangle_{q,t}$ we can write
\begin{equation} \label{4.1}
\prod_{i=1}^{n}\frac{ 1}{ (x_i ;q)_{r}}=
\sum_{\eta}\frac{\langle \prod_{i=1}^{n}\frac{ 1}{ (x_i ;q)_{r}},
E_{\eta}(x;q,t) \rangle_{q,t} }{
\langle E_{\eta}(x;q,t),E_{\eta}(x;q,t) \rangle_{q,t}}
E_{\eta}(x;q,t)
\end{equation}
Comparing (\ref{4.1}) with Proposition \ref{Bi-E} we have
\begin{equation}
{\rm CT}\left(  \prod_{i=1}^{n}\frac{ 1}{ (x_i ;q)_{r}}E_{\eta}(x^{-1};q^{-1},t^{-1}) 
W(x) \right)
=\frac{[q^{r}]_{\eta^+}^{q,t} }{ u_{\eta}(q,t)d_{\eta}(q,t)}
{\cal N}_{\eta}^{(E)}(q,t)
\end{equation}
Letting $x_{i}\mapsto x_{n-i+1}^{-1}$ inside the argument of
the constant term function, an operation that leaves it's
value unchanged, we obtain
\begin{equation} \label{4.2}
{\rm CT} \left(  \prod_{i=1}^{n}
\frac{ 1}{ (x_{i}^{-1} ;q)_{r}}{\hat{E}}_{\eta}(x;q,t) 
W(x) \right)
=\frac{[q^{r}]_{\eta^+}^{q,t} }{ u_{\eta}(q,t)d_{\eta}(q,t)}
{\cal N}_{\eta}^{(E)}(q,t)
\end{equation}
Our first task is to manipulate (\ref{4.2}) so that
$\prod_{i=1}^{n}(x_{i}^{-1} ;q)_{r}^{-1}$ is replaced by 
$\prod_{i=1}^{n}(x_{i};q)_{a}(qx_{i}^{-1} ;q)_{b}$.
We require
\begin{lemma} \label{xE}
We have  
\begin{equation} \label{xE1}
x^{p} E_{\eta}(x;q,t)= E_{\eta +p}(x;q,t)
\end{equation}
where $\eta +p=(\eta_{1} +p,\ldots ,\eta_{n} +p)$ and 
$x^{p}=(x_{1}\ldots x_{n})^{p}$
\end{lemma}
\noindent {\it Proof.} \quad From the definition of 
$Y_{i}$ we have
\begin{equation}
Y_{i} \, \left(x^{p}E_{\eta}(x;q,t)\right)= q^{p}x^{p}Y_{i}E_{\eta}(x;q,t)
\end{equation}
Using (\ref{E2}) we  obtain
\begin{eqnarray}
Y_{i} \, \left(x^{p}E_{\eta}(x;q,t)\right)
&=&
q^{p}t^{ {\bar{\eta}}_{i} }x^{p}E_{\eta}(x;q,t)
\\
&=&
t^{(\overline{\eta_{i}+p})_{i}}x^{p}E_{\eta}(x;q,t)
\end{eqnarray}
From the defining properties (\ref{E1}), (\ref{E2}) 
we then have the required conclusion.
\hfill $\Box$ \medskip

\begin{cor} \label{xdualE}
\begin{equation}
x^{p}{\hat{E}}_{\eta}(x;q,t) ={\hat{E}}_{\eta +p}(x;q,t) 
\end{equation}
\end{cor}

Using the above proof we can extend the non-symmetric 
Macdonald polynomials to include Laurent polynomials.
The defining properties of these Laurent polynomials $E_{\eta}$
are the same as for the ordinary non-symmetric 
Macdonald polynomials except that they are indexed
by compositions which can have negative parts. The
non-symmetric Macdonald Laurent polynomials can be
expressed in terms of the ordinary non-symmetric 
Macdonald polynomials by use of (\ref{xE1}).
The dual non-symmetric Macdonald polynomials can be
similarly extended to include Laurent polynomials.

Consider (\ref{4.2}) with $\eta$ replaced by $\eta +a$.
Using Lemma (\ref{xdualE}) we can write
${\hat{E}}_{\eta +a}=x^{a}{\hat{E}}_{\eta} $.
Set $r=-a-b$ with $-a,r \in {\mathbb{Z}}_{\leq 0}$.
A brief calculation shows that
\begin{equation}
x^{a}\prod_{i=1}^{n}\frac{1 }{ (x_{i}^{-1} ;q)_{r} }=
(-1)^{na}q^{-\frac{na }{ 2}(2b +a +1)}
\prod_{i=1}^{n}({x'}_{i};q)_{a}(\frac{q }{{x'}_{i}} ;q)_{b}
\end{equation} 
where ${x'}_{i}=q^{b+1}x_i$.
Substituting into (\ref{4.2}) we obtain
\begin{equation}
{\rm CT}\left(  \prod_{i=1}^{n}
({x}_{i};q)_{a}(\frac{q}{{x}_{i}} ;q)_{b}
{\hat{E}}_{\eta}(q^{-(b+1)}x;q,t) 
W(x) \right) 
=(-1)^{na}
\frac{q^{\frac{na }{ 2}(2b+a +1)}\, [q^{r}]_{\eta^{+} +a}^{q,t} 
}{ u_{\eta  +a}(q,t) \, d_{\eta +a}(q,t)}
{\cal N}_{\eta +a}^{(E)}(q,t) 
\end{equation}
Since ${\hat{E}}_{\eta}(cx)=c^{|\eta|}{\hat{E}}_{\eta}(x)$ and
${\cal N}_{\eta +a}^{(E)}(q,t)={\cal N}_{\eta}^{(E)}(q,t)$
we get
\begin{equation}
{\rm CT}\left(  \prod_{i=1}^{n}
({x}_{i};q)_{a}(\frac{q}{{x}_{i}} ;q)_{b}
{\hat{E}}_{\eta}(x;q,t) 
W(x) \right) =
(-1)^{na}\frac{ q^{\left( \frac{na }{ 2}(2b +a +1) + (b+1)|\eta| \right) }
\, [q^{r}]_{\eta^{+} +a}^{q,t} }{ u_{\eta  +a}(q,t) \, d_{\eta +a}(q,t)}
 \, {\cal N}_{\eta}^{(E)}(q,t)     \label{4.4}
\end{equation}
The dependence on $a$ in $1/ u_{\eta  +a} d_{\eta +a}$
can be determined using
\begin{lemma} \label{E_}
 We have  
\begin{equation}
 E_{\eta}(\frac{1 }{ x};q,t)= E_{- \underline{\eta}}(\underline{x};q,t)
\end{equation}
\end{lemma}
\noindent {\it Proof.} \quad
Let the star symbol ${\,}^{*}$ denote the involution on the ring of 
$n$-variable polynomials with coefficients
in ${\mathbb{Q}}(q,t)$ which sends $x_{i}\rightarrow  x_{i}^{-1}$,
$q \rightarrow q^{-1}$ and $t \rightarrow t^{-1}$.
Extend this operator to act on operators so that
for any operator $T$ and polynomial $f$, $T^{*} f^{*}= (Tf)^{*}$.
From the relations $T_{i}^{*}=T_{i}^{-1}$, $w^{*}=w$,
${\hat{T}}_{i}=T_{n-i}^{-1}$ and $\hat{w}=w^{-1}$ \cite{BF97h}
it follows that
\begin{equation} \label{4.3}
({Y_{n-i+1}^{*}})^{-1}=t^{1-n}\,{\hat{Y}}_{i}
\end{equation}
From (\ref{E2}) 
\begin{equation}
Y_{i}^{-1} \, E_{\eta}(x;q,t) = t^{-\bar{\eta}_i}\,E_{\eta}(x;q,t)
\end{equation}
Applying the ${\,}^{*}$ operator and replacing $i$ with
$n-i+1$ we get
\begin{equation}
(Y_{n-i+1}^{*})^{-1} \, E_{\eta}(x^{-1};q^{-1},t^{-1})
= t^{\bar{\eta}_{n-i+1}}\,E_{\eta}(x^{-1};q^{-1},t^{-1})
\end{equation}
Using (\ref{4.3}) we  obtain
\begin{equation}
{\hat{Y}}_{i} \, E_{\eta}(x^{-1};q^{-1},t^{-1})
= t^{n-1+\bar{\eta}_{n-i+1}}\,E_{\eta}(x^{-1};q^{-1},t^{-1})
\end{equation}
From the defining properties  (\ref{dual}) and (\ref{dual.b})
it follows that $E_{\eta}(x^{-1};q^{-1},t^{-1})$
is a dual non-symmetric Macdonald polynomial.
Since $E_{\eta}(x^{-1};q^{-1},t^{-1})$
has the same leading term as ${\hat{E}}_{-\underline{\eta}}(x;q,t)$
\begin{equation}
E_{\eta}(\frac{1 }{ x};q^{-1},t^{-1})=
E_{- \underline{\eta}}(\underline{x};q^{-1},t^{-1})
\end{equation}
The conclusion follows.
\hfill $\Box$ \medskip

\begin{cor}\label{E_dual}
We have
\begin{equation}
{\hat{E}}_{\eta}(\frac{1 }{ x};q,t)= {\hat{E}}_{- \underline{\eta}}(\underline{x};q,t)
\end{equation}
\end{cor}

Now
\begin{eqnarray} \nonumber
{\rm CT}\left(  \prod_{i=1}^{n}
({x}_{i};q)_{a}(\frac{q}{{x}_{i}} ;q)_{b}
{\hat{E}}_{\eta}(x;q,t) 
W(x) \right) 
&=&
{\rm CT}\left(  \prod_{i=1}^{n}
(\frac{q}{{x}_{i}} ;q)_{a}({x}_{i};q)_{b}
{\hat{E}}_{\eta}(\frac{q}{\underline{x}};q,t) 
W(\frac{q}{\underline{x}}) \right) 
\\
&=&
q^{|\eta|}{\rm CT}\left(  \prod_{i=1}^{n}
(\frac{q}{{x}_{i}} ;q)_{a}({x}_{i};q)_{b}
{\hat{E}}_{-\underline{\eta}}(x;q,t) 
W(x) \right) 
\end{eqnarray}
To obtain the first equality we have used the
invariance of the constant term identity
under $x_{i} \mapsto \frac{q}{{x}_{n-i+1}}$,
while to get the second equality we have used Corrollary \ref{E_dual}
and $W(\frac{q}{\underline{x}})=W(x)$.

Applying (\ref{4.4}) with $\eta$ replaced by $-\underline{\eta}$
and $a$ interchanged with $b$ gives
\begin{equation} \label{4.5}
{\rm CT}\left(  \prod_{i=1}^{n}
({x}_{i};q)_{a}(\frac{q}{{x}_{i}} ;q)_{b}
{\hat{E}}_{\eta}(x;q,t) 
W(x) \right) 
=
(-1)^{nb}\frac{q^{\left( \frac{nb }{ 2}(2a +b +1) -a|\eta| \right) } \,
[q^{r}]_{-{\underline{\eta}}^{+} +b}^{q,t}
 }{ u_{-{\underline{\eta}} +b}(q,t)\, 
 d_{-{\underline{\eta}} +b}(q,t)}
{\cal N}_{\eta}^{(E)}(q,t) 
\end{equation}
We write $\eta\leq c$ if $\eta_{i}\leq c$ for all $i=1,\ldots,n$. 
The equation (\ref{4.5}) is valid for $\eta \leq b$ 
while (\ref{4.4}) is valid for $\eta \geq -a$.
Equating the right hand sides of (\ref{4.4})
and (\ref{4.5}) and setting $a=0$
we obtain for $0 \leq \eta \leq b$
\begin{equation} \label{4.6}
\frac{1}{ u_{-{\underline{\eta}} +b}(q,t)\,
 d_{-{\underline{\eta}} +b}(q,t)}
=(-1)^{nb}
\frac{q^{\left( -\frac{nb }{ 2}(b+1)+ (b+1)|\eta| \right) }\,[q^{-b}]_{\eta+}
}{
u_{\eta}(q,t)\, d_{\eta}(q,t) \,[q^{-b}]_{-{\underline{\eta}}^{+} +b} }
\end{equation}
We can use (\ref{4.6}) to define $1/u_{-{\underline{\eta}} +b}
d_{-{\underline{\eta}} +b}$ for ${-{\underline{\eta}} +b} \leq 0$.
It then follows that (\ref{4.6}) is true for all 
$\eta \in {\mathbb{Z}}^{n}$
which in turn can be used to show that  (\ref{4.5}) 
is true for all $\eta \in {\mathbb{Z}}^{n}$.
Substituting (\ref{4.6}) into (\ref{4.5}) we then
obtain for $a,b \in {\mathbb{Z}}_{\geq 0}$
\begin{equation}  \label{4.gen}
{\rm CT}\left(  \prod_{i=1}^{n}
({x}_{i};q)_{a}(\frac{q}{{x}_{i}} ;q)_{b}
{\hat{E}}_{\eta}(x;q,t) 
W(x) \right) 
=\frac{  q^{\left( nab +(b+1-a)|\eta| \right) }
[q^{-b}]_{\eta+}\,[q^{-a-b}]_{-{\underline{\eta}}^{+} +b} }{
u_{\eta}(q,t)\, d_{\eta}(q,t) \,[q^{-b}]_{-{\underline{\eta}}^{+} +b}}
{\cal N}_{\eta}^{(E)}(q,t) 
\end{equation}
To extend this result to all $a,b \in {\mathbb{R}}$ we note that
both sides of (\ref{4.gen}) are series in $q,q^{a},q^{b}$ and 
$q^{\lambda}$. We then apply Lemma\ \ref{Stem} twice, once with respect to
$q^{a}$ and once with respect to $q^{b}$.

The identity (\ref{4.gen}) can be simplified
by taking the limit
$a \rightarrow \infty$ with $r=-a-b$ remaining constant

For this purpose it is
convenient to first take the ratio of (\ref{4.gen})
to that obtained with $\eta =0$, thus obtaining
\begin{equation}  \label{4.gen1}
 \frac{{\rm CT} \left(  \prod_{i=1}^{n}
({x}_{i};q)_{a}    (\frac{q}{{x}_{i}} ;q)_{b}
{\hat{E}}_{\eta}(x;q,t) 
W(x) \right)
}{
{\rm CT} \left(  \prod_{i=1}^{n}
({x}_{i};q)_{a} (\frac{q}{{x}_{i}} ;q)_{b}
W(x) \right)}
=
\frac{q^{(b+1-a)|\eta|}\,[q^{-b}]_{\eta+}\, [q^{-b}]_{b}
\, [q^{-a-b}]_{-{\underline{\eta}}^{+} +b}\, {\cal N}_{\eta}^{(E)}(q,t) 
}{
u_{\eta}(q,t)\, d_{\eta}(q,t) \,[q^{-a-b}]_{b}
[q^{-b}]_{-{\underline{\eta}}^{+} +b}\,{\cal N}_{0}^{(E)}(q,t) }
\end{equation}
where we used the facts that $[q^{-b}]_{0}=d_{0}=u_{0}=1$.
Computing the the asymptotics requires 
\begin{lemma} {\rm\cite{BF98a}} \label{q-int}
For a general Laurent polynomial $f(s_{1}, \ldots ,s_{n})$
we have
$$
\left( \frac{\Gamma_{q}[a+1]}{ \Gamma_{q}[-b] \, \Gamma_{q}[a+b+1] }
 \right)^{n}\prod_{i=1}^{n}\int_{0}^{1}d_{q}s_{i}s_{i}^{-b-1}
\frac{(qs_{i};q)_{\infty} }{ (q^{a+b+1}s_{i};q)_{\infty}}
f(s_{1}, \ldots ,s_{n})  \qquad
$$
\begin{equation}  \label{q-int.1}
\qquad\qquad\qquad\qquad =
\left( \frac{(q,q)_{a}\, (q,q)_{b}}{(q,q)_{a+b} }
 \right)^{n}
{\rm CT}_{\{ s \}}\left(\prod_{i=1}^{n}(s_{i},q)_{a}\, 
(\frac{q }{ s_{i}},q)_{b}
f(q^{-(b+1)}s_{1}, \ldots ,q^{-(b+1)}s_{n})\right)
\end{equation}
where 
$\int_{0}^{1}f(s)d_{q}s:=(1-q)\sum_{i=0}^{\infty}f(q^{j})q^{j}$
is the $q$-integral.
\end{lemma}

\noindent {\it Remark.} \quad There is a typing error in
the statement of the above lemma in {\rm\cite{BF98a}}.
\medskip
 
\noindent Lemma \ref{q-int} allows us to deduce
\begin{lemma}  \label{assym}
Letting $\lambda \in Z_{\geq 0}$ and $a+b=\mbox{const}$
we have
\begin{equation} \label{assym.1}
\lim_{a \rightarrow \infty}
\frac{{\rm CT}\left(  \prod_{i=1}^{n}
({x}_{i};q)_{a}(\frac{q}{{x}_{i}} ;q)_{b}
{\hat{E}}_{\eta}(x;q,t) 
W(x)\right)
}{
{\rm CT}\left(  \prod_{i=1}^{n}
({x}_{i};q)_{a}     (\frac{q }{{x}_{i} } ;q)_{b}
W(x) \right)}
=
q^{|\eta|(b+1)}{\hat{E}}_{\eta}(t^{\underline{\delta}};q,t) 
\end{equation}
\end{lemma}
\noindent {\it Proof.} \quad 
Fixing $r=-a-b $ and applying Lemma (\ref{q-int})
we obtain
\begin{equation} \label{assym.0}
\frac{{\rm CT}\left(  \prod_{i=1}^{n}
({x}_{i};q)_{a}(\frac{q}{{x}_{i}} ;q)_{b}
{\hat{E}}_{\eta}(x;q,t) 
W(x) \right)
}{
{\rm CT}\left(  \prod_{i=1}^{n}
({x}_{i};q)_{a}(\frac{q}{{x}_{i}} ;q)_{b}
W(x) \right)}
=q^{(b+1)|\eta|}
\frac{ \prod_{i=1}^{n}\int_{0}^{1}d_{q}s_{i}s_{i}^{-b-1}
\frac{ (qs_{i};q)_{\infty} }{ (q^{-r+1}s_{i};q)_{\infty} }
{\hat{E}}_{\eta}(s;q,t) 
W(s)
}{
 \prod_{i=1}^{n}\int_{0}^{1}d_{q}s_{i}s_{i}^{-b-1}
\frac{ (qs_{i};q)_{\infty}}{ (q^{-r+1}s_{i};q)_{\infty} }
W(s) }
\end{equation}
Using the definition of the $q$-integral we have
\begin{equation}
\prod_{i=1}^{n}\int_{0}^{1}d_{q}s_{i}s_{i}^{-b-1}
\frac{ (qs_{i};q)_{\infty} }{ (q^{-r+1}s_{i};q)_{\infty} }
{\hat{E}}_{\eta}(s;q,t) 
W(s)  \qquad\qquad\qquad\qquad\qquad\qquad\qquad
\end{equation}
\begin{equation} \qquad\qquad\qquad
=(1-q)^{n}\sum_{k_{i} \in {\mathbb{Z}}_{\geq 0}}
q^{ -b\sum_{i=1}^{n}k_{i} }
\prod_{i=1}^{n}
\frac{ (q^{ k_{i}+1 };q)_{\infty} }{ (q^{ (2-r)k_{i} };q)_{\infty} }
W(q^{k_{1}},\ldots,q^{k_{n}}) {\hat{E}}_{\eta}(q^{k_{1}},\ldots,q^{k_{n}};q,t)
\end{equation}

Suppose $\lambda \in {\mathbb{Z}}_{\geq 0}$ .
Then 
\begin{equation}
W(q^{k_{1}},\ldots,q^{k_{n}})=0 \quad \mbox{if} \quad
k_{i+1}=k_{i}-\lambda, \ldots,k_{i}+\lambda - 1 \quad
\mbox{while}
\quad
W(1,q^{\lambda},\ldots,q^{\lambda (n-1)})\neq 0
\end{equation}
It follows that in the limit $b \rightarrow \infty$ with 
$r=-a-b$ fixed
\begin{displaymath}
\prod_{i=1}^{n}\int_{0}^{1}d_{q}s_{i}s_{i}^{-b-1}
\frac{ (qs_{i};q)_{\infty} }{ (q^{-r+1}s_{i};q)_{\infty} }
{\hat{E}}_{\eta}(s;q,t) 
W(s)   \qquad \qquad \qquad \qquad \qquad \qquad
\qquad \qquad
\end{displaymath}
\begin{equation}  \label{assym2}
\qquad \qquad  \sim (1-q)^{n}
q^{ -b\sum_{i=1}^{n}\lambda (i-1) }
\prod_{i=1}^{n}
\frac{ (q^{ \lambda (i-1)+1 };q)_{\infty} 
}{ (q^{ (2-r)\lambda (i-1) };q)_{\infty} }
{\hat{E}}_{\eta}(t^{\underline{\delta}};q,t)
W(t^{\underline{\delta}})
\end{equation}
Substituting (\ref{assym2}) in to (\ref{assym.0}) gives (\ref{assym.1}).
\hfill $\Box$ \medskip

Lemma \ref{assym} gives the asymptotics of the left hand
side of (\ref{4.gen1}).
We now seek the asymptotics of the right hand side 
of (\ref{4.gen1}).
It follows from the property
\begin{equation}  \label{gamma1}
\Gamma_{q}[x+1]=[x]_{q}\,\Gamma_{q}[x]  \qquad \mbox{where} 
\quad [x]_{q}:=\frac{1-q^{x} }{1-q}
\end{equation}
that for all $k \in {\mathbb{Z}}$
\begin{equation}  \label{gamma2}
\frac{\Gamma_{q}[x+k] }{ \Gamma_{q}[x]}
=(-1)^{k}q^{kx +\frac{1 }{ 2}k(k-1)}
\frac{\Gamma_{q}[1-x] }{ \Gamma_{q}[1-(x+k)]}
\end{equation}

 Using these properties  
along
with 
\begin{equation}
\frac{\Gamma_{q}[x+a] }{ \Gamma_{q}[x]} \sim [x]^{a}_{q}
\qquad \mbox{as} \quad x \rightarrow \infty
\end{equation}
shows that in the limit $a \rightarrow \infty$ with $a+b$ fixed
\begin{eqnarray}
[q^{-b}]_{\eta^{+}}^{q,t} 
&\sim &
t^{l'(\eta^{+})}(1-q)^{|\eta|}
[-b]_{q}^{|\eta|}
\\
\frac{[q^{-a-b}]_{-{\underline{\eta}}^{+}+b}^{q,t}
}{
[q^{-a-b}]_{b}^{q,t}} 
&\sim&
(-1)^{|\eta|}(1-q)^{-|\eta|}
q^{\left( a|\eta|+ \frac{1 }{ 2}
\sum_{i=1}^{n}\eta_{i}^{+}(\eta_{i}^{+}+1)\right)}
[a]_{q}^{-|\eta|}
\\  \label{4.simp}
\frac{[q^{-b}]_{b}^{q,t} }{ [q^{-b}]_{-{\underline{\eta}}^{+}+b}^{q,t} }
&=&
(-1)^{|\eta|}
q^{ -\frac{1 }{ 2}
\sum_{i=1}^{n}\eta_{i}^{+}(\eta_{i}^{+}+1)}
t^{-l'(\eta^{+})}[q^{1+\lambda (n-1)}]_{\eta^{+}}^{q,t}
\end{eqnarray}
Substituting these results in to the right hand side of
(\ref{4.gen1}) and using Lemma \ref{assym} we have in the limit
$a \rightarrow \infty$ with $a+b$ fixed
and $\lambda \in {\mathbb{Z}}_{\geq 0}$
\begin{equation}  \label{Eu}
{\hat{E}}_{\eta}(t^{\underline{\delta}};q,t)
=
\frac{[q^{1+\lambda (n-1)}]_{\eta^{+}}^{q,t}
}{
u_{\eta}(q,t) \, d_{\eta}(q,t)}
\frac{{\cal N}_{\eta}^{(E)}(q,t) 
}{  {\cal N}_{0}^{(E)}(q,t) }
\end{equation}
Since both sides of this expression can be written as power
series in $q$ and $t$ for  $0 <q,\,t<1$ we can apply
Lemma \ref{Stem} to extend 
the validity of this result to all $\lambda >0$.
Using this result,  (\ref{4.simp}) and
\begin{equation}
\frac{[q^{-a-b}]_{-{\underline{\eta}}^{+}+b}^{q,t}
}{
[q^{-a-b}]_{b}^{q,t}}
=\frac{(-1)^{|\eta|}t^{l'(\eta^+)}
q^{\left( a|\eta| + \frac{1 }{ 2} 
\sum_{i=1}^{n}\eta_{i}(\eta_{i}+1) \right) }
}{
[q^{1+a+ \lambda (n-1)}]_{\eta^{+}} }
\end{equation}
we can simplify (\ref{4.gen1}) to obtain
\begin{equation} \label{genmorris.-1}
\frac{{\rm CT}\left(  \prod_{i=1}^{n}
({x}_{i};q)_{a}(\frac{q}{{x}_{i}} ;q)_{b}
{\hat{E}}_{\eta}(x;q,t) 
W(x) \right)
}{
{\rm CT}\left(  \prod_{i=1}^{n}
({x}_{i};q)_{a}(\frac{q}{{x}_{i}} ;q)_{b}
W(x) \right)}
= 
q^{(b+1)|\eta|}{\hat{E}}_{\eta}(t^{\underline{\delta}};q,t)
\frac{[q^{-b}]_{\eta^{+}}   }{
[q^{1+a+ \lambda (n-1)}]_{\eta^{+}} }
\end{equation} 
It follows from Lemma \ref{E-T-hatE} and 
(\ref{TE}) that 
 $\{ E_{\eta} \}_{\eta^+ =\ka}$ 
and $\{ {\hat{E}}_{\eta} \}_{\eta^+ =\ka}$ 
span the same set of functions.
In particular we can write
\begin{equation}
E_{\mu}(x;q,t) =\sum_{\{\eta|\eta^+ = \mu^+\}}
c_{\mu\eta}{\hat{E}}_{\eta}(x;q,t) 
\end{equation}
for scalars $c_{\mu\eta}$.
Multiplying both sides of (\ref{genmorris.-1}) by $c_{\mu\eta}$
and summing over distinct permutations of $\mu^+$ we obtain

\begin{prop}  \label{genmorris}
\begin{equation}  \label{genmorris1}
\frac{{\rm CT}\left(  \prod_{i=1}^{n}
({x}_{i};q)_{a}(\frac{q}{{x}_{i}} ;q)_{b}
{{E}}_{\eta}(x;q,t) 
W(x) \right)
}{
{\rm CT}\left(  \prod_{i=1}^{n}
({x}_{i};q)_{a}(\frac{q}{{x}_{i}} ;q)_{b}
W(x) \right)}
= 
q^{(b+1)|\eta|}{{E}}_{\eta}(t^{\underline{\delta}};q,t)
\frac{[q^{-b}]_{\eta^{+}}
}{
[q^{1+a+ \lambda (n-1)}]_{\eta^{+}}}
\end{equation}   
\end{prop}
Note that by multiplying both sides of  (\ref{genmorris1})
by ${d'}_{\eta^{+}}(q,t)/{d'}_{\eta}(q,t)$, summing over distinct
permutations of $\ka = \eta^+$ and applying (\ref{P-E})
we get back Proposition \ref{genmorris} with ${{E}}_{\eta}$
replaced by the symmetric Macdonald polynomial $P_{\ka}$. 
Restraining $\lambda$ to be a non-negative
integer we can use
 Lemma \ref{q-int} to transform (\ref{genmorris1})
into a generalisation of the $q$-Selberg integral.

\begin{prop} \label{genqselprop1}
\begin{displaymath}
\prod_{i=1}^{n}\int_{0}^{1}d_{q}s_{i}s_{i}^{x-1}
\frac{(qs_{i};q)_{\infty} }{ (q^{a+b+1}s_{i};q)_{\infty}}
{{E}}_{\eta}(s;q,t)
\prod_{i<j}s_{i}^{2\lambda}(q^{1-\lambda}
\frac{s_{j}}{ s_{i}};q)_{2\lambda}
\qquad\qquad\qquad\qquad
\qquad\qquad\qquad\qquad
\end{displaymath}
\begin{equation}  \label{genqsel1}
\qquad\qquad
={{E}}_{\eta}(t^{\underline{\delta}};q,t)
\frac{[q^{x+\lambda (n-1)}]_{\eta^{+}}
}{
[q^{x+y+ 2\lambda (n-1)}]_{\eta^{+}}}
\prod_{i=1}^{n}\int_{0}^{1}d_{q}s_{i}s_{i}^{x-1}
\frac{(qs_{i};q)_{\infty} }{ (q^{a+b+1}s_{i};q)_{\infty}}
\prod_{i<j}s_{i}^{2\lambda}(q^{1-\lambda}\frac{s_{j}}{ s_{i}};q)_{2\lambda}
\end{equation}
\end{prop}

\noindent {\it Proof.} \quad
Apply 
(\ref{q-int.1}) to (\ref{genmorris1}) and
write 
\begin{equation}
W(s)=(-1)^{\lambda n(n-1)/2}q^{n(n-1)\lambda (\lambda -1)/4}
\left(
\prod_{i=1}^{n}s_{i}^{-\lambda (n-1)}\right)
\prod_{i<j}s_{i}^{2\lambda}(q^{1-\lambda}\frac{s_{j}}{ s_{i}};q)_{2\lambda}
\end{equation}
Then let $x=-b-\lambda (n-1)$, $y=a+b+1$.
\hfill $\Box$ \medskip

The above derivation of Proposition \ref{genqselprop1}
has some further consequences in relation to the general theory.
First, it allows new derivations of the $q$-Morris identity
and the $q$-Selberg integral.

\begin{prop} {\rm \cite{Mo82}}\label{q-Morris-prop}
\begin{equation}  \label{q-morris}
{\rm CT}\left(  \prod_{i=1}^{n}
({x}_{i};q)_{a}(\frac{q}{{x}_{i}} ;q)_{b}
W(x) \right) 
=\frac{\Gamma_{q}[1+a+b+\lambda (i-1)] \Gamma_{q}[1+ \lambda i] 
}{
 \Gamma_{q}[1+a+\lambda (i-1)] \Gamma_{q}[1+b+\lambda (i-1)]
\Gamma_{q}[1+\lambda]}
\end{equation} 
\end{prop}

\noindent {\it Proof. }
Letting $\eta=0$ in (\ref{4.gen}) and using (\ref{cajamarca}) we have
\begin{equation}  
{\rm CT}\left(  \prod_{i=1}^{n}
({x}_{i};q)_{a}(\frac{q}{{x}_{i}} ;q)_{b}
W(x) \right) 
=
q^{abn}{\cal N}_{0}^{(E)}(q,t)
\prod_{i=1}^{n}\frac{\Gamma_{q}[-a-\lambda (i-1)]\, 
\Gamma_{q}[-b-\lambda (i-1)]
}{
\Gamma_{q}[-a-b-\lambda (i-1)]\, \Gamma_{q}[-\lambda (i-1)]}
\end{equation}
The $q$-Morris identity (\ref{q-morris}) is then obtained by 
using the properties
(\ref{gamma1}), (\ref{gamma2}) and the evaluation \cite{Ad77}
\begin{equation}
{\cal N}_{0}^{(E)}(q,t) 
=\frac{\Gamma_{q}[\lambda n +1] }{ \Gamma_{q}[\lambda +1]^{n}}
\end{equation}
\hfill $\Box$ \medskip

The $q$-Selberg integral can be evaluated
as in \cite{BF98a} by applying Lemma \ref{q-int} 
to the $q$-Morris identity
and making some manipulations.

\begin{prop} {\rm\cite{As80}}
\begin{eqnarray} \nonumber
&&\prod_{i=1}^{n}\int_{0}^{1}d_{q}s_{i}s_{i}^{x-1}
\frac{(qs_{i};q)_{\infty} }{ (q^{a+b+1}s_{i};q)_{\infty}}
\prod_{i<j}s_{i}^{2\lambda}(q^{1-\lambda}\frac{s_{j}}{ s_{i}};q)_{2\lambda}
\qquad\qquad\qquad\qquad\qquad
\\  \label{qsel} &&\qquad\qquad\qquad
=q^{
\lambda x \binom{n}{2}
+2{\lambda}^{2}\binom{n}{3}}
\prod_{i=1}^{n}\frac{\Gamma_{q}[x+\lambda (i-1)]\Gamma_{q}
[y+\lambda (i-1)]
\Gamma_{q}[1+\lambda i]}{ \Gamma_{q}[x+y+\lambda (n+i-2)]
\Gamma_{q}[\lambda +1]}
\end{eqnarray}
\end{prop}
 
We can use (\ref{qsel}) to simplify Proposition \ref{genqselprop1}.

\begin{prop}
\begin{eqnarray} \nonumber
&&\prod_{i=1}^{n}\int_{0}^{1}d_{q}s_{i}s_{i}^{x-1}
\frac{(qs_{i};q)_{\infty} }{ (q^{a+b+1}s_{i};q)_{\infty}}
E_{\eta}(s;q,t)
\prod_{i<j}s_{i}^{2\lambda}
(q^{1-\lambda}\frac{s_{j}}{ s_{i}};q)_{2\lambda}
\qquad\qquad\qquad\qquad
\qquad\qquad\qquad\qquad\qquad
\\ \label{genqsel}
&&\qquad\qquad\qquad
=q^{\lambda x \binom{n}{2}
+2{\lambda}^{2}\binom{n}{3}}
E_{\eta}(t^{\bar{\delta}};q,t)
\prod_{i=1}^{n}\frac{\Gamma_{q}[\lambda i+1]\,
\Gamma_{q}[x+\lambda (n-i)+\eta^{+}_{i}]\,
\Gamma_{q}[y+\lambda (n-i)]
}{ \Gamma_{q}[\lambda +1]\,
\Gamma_{q}[x+y+\lambda (2n-i-1)+ \eta^{+}_{i}]}
\end{eqnarray} 
\end{prop}
This formula is
a generalisation of the integration
formula of Kadell \cite{Ka88} and Kaneko \cite{Kan96}. 
The formula of \cite{Kan96} can be reclaimed by
multiplying both sides of (\ref{genqsel}) by 
${d'}_{\eta^{+}}(q,t)/{d'}_{\eta}(q,t)$
and summing over distinct
permutations of $\ka = \eta^+$ 
using (\ref{dsymsum}).

The second consequence of the derivation of
Proposition \ref{genqselprop1} is that it allows
us to calculate the normalisation constant $u_{\eta}(q,t)$
appearing in (\ref{E-sum1}).

\begin{prop} {\rm\cite{MN}} \label{u-eta}
\begin{equation} \label{u-eta.1}
u_{\eta}(q,t)=
\frac{{d'}_{\eta}(q,t) }{ {d}_{\eta}(q,t) }
\end{equation}
\end{prop}

\noindent {\it Proof.}
Using (\ref{e-G}),  (\ref{E(t)}), (\ref{e/d(1/q,1/t)}) and 
 $E_{\eta}(cx)=c^{|\eta|}E_{\eta}(x)$
we obtain 
$ u_{\eta}(q,t)=
{d'}_{\eta}(q,t) / {d}_{\eta}(q,t) $ for 
$\lambda \in {\mathbb{Z}}_{\geq 0}$.
Since both sides of this expression can be written
as formal power series in $q$ and $t$ if $0 <q,\,t<1$,
we can apply Lemma \ref{Stem} to show that this result is 
true for all $\lambda > 0$.
\hfill $\Box$ \medskip

\setcounter{equation}{0}
\section{\mbox{Symmetric Macdonald Polynomial Theory}}

In this section we will deduce analogues of Propositions
$3.1-3.5$ and Proposition \ref{u-eta} for the symmetric Macdonald
polynomials. This will be done by exploiting the relationships
between the symmetric, $q$-antisymmetric and non-symmetric
Macdonald polynomials.

In order to deduce the analogue of Proposition \ref{E-sum}
we need to derive the following two results.
The first reveals the relationship between the symmetric
and $q$-antisymmetric Macdonald polynomials.

\begin{lemma} 
\label{s-qa}
\begin{equation}  \label{s-qa1}
S_{\ka + \de}(x;q,t)=t^{-\frac{n(n-1) }{ 2}}\Delta_{t}(x)P_{\ka}(x;q,qt)
\end{equation}
\end{lemma}

\noindent {\it Proof.}\quad
Kadell's Lemma \cite{Ka88} gives for any antisymmetric function $h(x)$
\begin{equation}
{\rm CT}\left( \prod_{i < j}(x_i -ax_j)h(x)\right)=\frac{[n]_{a}! }{ n!}
{\rm CT}\left( \prod_{i < j}(x_i -x_j)h(x)\right)
\end{equation}
Consider 
\begin{equation} \label{5.1a}
\langle \Delta_{t}(x)P_{\ka}(x;q,qt),
\Delta_{t}(x)P_{\lambda}(x;q,qt) \rangle_{q,t}
={\rm CT}\left( \prod_{i < j}(x_i -\frac{1 }{ t}x_j)h(x)\right)
\end{equation}
where
\begin{equation} 
h(x):=\prod_{i < j}(\frac{1 }{ x_{i}}-\frac{1 }{ x_{j}})
(q\frac{x_{i} }{ x_{j}};q)_{\lambda}
(q\frac{x_{j} }{ x_{i}};q)_{\lambda}
P_{\ka}(x;q,qt)P_{\lambda}
(\frac{1 }{ x};\frac{1 }{ q},\frac{1 }{ qt })
\end{equation}
is an antisymmetric polynomial.

Applying Kadell's Lemma twice to the left hand side of
(\ref{5.1a}) gives
\begin{eqnarray}
\langle\Delta_{t}(x)P_{\ka}(x;q,qt),
\Delta_{t}(x)P_{\lambda}(x;q,qt) \rangle_{q,t}
&=&
\frac{[n]_{t^{-1}}! }{ [n]_{tq}!}
{\rm CT}\left( \prod_{i < j}(x_i -qtx_j)h(x)\right)
\\
&=&
\frac{[n]_{t^{-1}}! }{ [n]_{tq}!}
\langle P_{\ka}(x;q,qt),
P_{\lambda}(x;q,qt) \rangle_{q,qt}
\\
&=&
\frac{[n]_{t^{-1}}! }{ [n]_{tq}!}
\langle P_{\ka}(x;q,qt),
P_{\ka}(x;q,qt) \rangle_{q,qt}\,\de_{\ka \lambda}
\end{eqnarray}
The polynomials 
$ t^{-n(n-1) /2}\Delta_{t}(x)P_{\ka}(x;q,qt) $ 
then form an 
orthogonal set with respect to $\langle \cdot, \cdot \rangle_{q,t}$.
Since they also satisfy (\ref{S2}) 
with leading term ${m'}_{\ka +\delta}$
we obtain (\ref{s-qa1}).
\hfill $\Box$ \medskip

For the second result define an 
equivalence relationship $\acute{\sim}$ such that
\begin{equation}
f(x)\, \acute{\sim} \, g(x) \quad \mbox{iff} \quad f(x)-g(x)=
\sum_{i}x^{\eta^{(i)}}  \quad
\mbox{where all the $\eta^{(i)}$ have repeated parts.}
\end{equation}
Note that if $f(x)\, \acute{\sim} \, g(x)$ then it follows from
(\ref{U=0}) that $U^{-}f(x)=U^{-}g(x)$.
The sought identity is a partial confirmation of a $q$-generalisation
of the Cauchy double alternant formula.
\begin{lemma}  \label{Cauchy}
\begin{equation}  \label{Cauchy1}
{U^{-}}^{(x)}\left( \prod_{i=1}^{n}\frac{1}{ 1-tx_{i}y_{i}}
\prod_{j<i}^{n}\frac{ 1-x_{i}y_{j}}{  1-tx_{i}y_{j}}
\right)
=
\frac{F(y) \Delta_{t}(x) }{\prod_{i,j}^{n}(1-tx_{i}y_{j})}
\end{equation}
where $F(y)\, \acute{\sim} \,\Delta_{t^{-1}}(y)$.
\end{lemma}

\noindent {\it Remark.}
We shall see later (\ref{2.Cauchy.a}) that $F(y)=\Delta_{t^{-1}}(y)$.
\medskip

\noindent 
{\it Proof.} \quad
We shall first show that
\begin{equation} \label{Cauchy1.5}
\Delta_{t}(y)\, \acute{\sim} \,(-1)^{n(n-1)/2}\sum_{\sigma \in S_{n}}
(-t)^{l(\sigma)}y^{\sigma -1}
\end{equation}
where for the permutation $\sigma=(\sigma (1), \cdots,\sigma (n))$ 
\begin{equation}
y^{\sigma -1}:=y_{1}^{\sigma(1)-1}\ldots y_{n}^{\sigma(n)-1}
=y_{\sigma^{-1}(n)}^{n-1}\ldots y_{\sigma^{-1}(1)}^{0}
\end{equation}
It is clear that the only terms of $\Delta_{t}(y)$ with
non-repeating parts are $\{y^{\sigma -1} \}$.
\begin{equation}
\Delta_{t}(y)
\, \acute{\sim} \,\sum_{\sigma \in S_{n}}a_{\sigma}y^{\sigma -1}
\end{equation}
Also, given a permutation $\sigma$, we can write
\begin{equation}
\prod_{i<j}(ty_{i}-x_{j})\, \acute{\sim} \,
(-1)^{n(n-1)/2}
\prod_{i=1}^{n}\left(\prod_{
\stackrel{k<\sigma^{-1}(i)}{k\neq \sigma^{-1}(n),\ldots,\sigma^{-1}(i+1)}}
(y_{\sigma^{-1}(i)}-ty_{k}) \right)
\left( \prod_{
\stackrel{k>\sigma^{-1}(i)}{k\neq \sigma^{-1}(n),\ldots,\sigma^{-1}(i+1)}}
(y_{k}-ty_{\sigma^{-1}(i)}) \right)
\end{equation}
It is then apparent that the coefficient of
$y^{\sigma} -1$ is
\begin{eqnarray}  \nonumber
a_{\sigma}&=&(-1)^{n(n-1)/2}\prod_{i=1}^{n}
(-t)^{\left
(n-\sigma^{-1}(i)-\#\{j:\sigma^{-1}(i)< \sigma^{-1}(j),i<j\}
\right)}
\\  \nonumber
&=&(-1)^{n(n-1)/2}(-t)^{n(n-1)/2-
\#\{(i,j):\sigma(i)< \sigma(j),i<j\}}
\\
&=&(-1)^{n(n-1)/2}(-t)^{l(\sigma)}
\end{eqnarray}
so (\ref{Cauchy1.5}) follows.

Let us now consider the left hand side of (\ref{Cauchy1}).
We can write
\begin{equation}
\prod_{i=1}^{n}\frac{1}{ 1-tx_{i}y_{i}}
\prod_{j<i}^{n}\frac{ 1-x_{i}y_{j}}{  1-tx_{i}y_{j}}
=
\left( \prod_{i,j=1}^{n}\frac{1 }{ 1-tx_{i}y_{j}}\right)
\prod_{i \neq j}(1-t^{\theta (j-i)}x_{i}y_{j})
\end{equation}
where
\begin{equation}
\theta(s)=
\left\{
\begin{array}{ll}
1 & s>0
\\
0 & s<0
\end{array}
\right.
\end{equation}
Since symmetric functions commute with $U^{-}$the left hand side
of (\ref{Cauchy1}) can be written 
\begin{equation}  \label{Cauchy3}
\left( \prod_{i,j=1}^{n}\frac{1 }{ 1-tx_{i}y_{j}}\right)
{U^{-}}^{(x)}\prod_{i \neq j}(1-t^{\theta (j-i)}x_{i}y_{j})
\end{equation}
For the power series $f=\sum_{\eta}c_{\eta}x^{\eta}$
let $[x^{\eta}]f $ denote $c_{\eta}$. Since the only terms 
 of $\prod_{i \neq j}(1-t^{\theta (j-i)}x_{i}y_{j})$
with $x^{\eta}$ having non-repeating parts are $\{ x^{\sigma -1} \}$
we have
\begin{eqnarray}  \nonumber
{U^{-}}^{(x)}\prod_{i \neq j}(1-t^{\theta (j-i)}x_{i}y_{j})
&=& 
\sum_{\sigma \in S_{n}}\left( [ x^{\sigma -1}]\,
\prod_{i \neq j}(1-t^{\theta (j-i)}x_{i}y_{j}) \right)
{U^{-}}^{(x)} x^{\sigma -1}
\\ &=&  \label{Cauchy4}
(-t)^{-n(n-1)/2}\Delta_{t}(x)
\sum_{\sigma \in S_{n}}(-1)^{-l(\sigma)}\,[ x^{\sigma -1}]\,
\prod_{i \neq j}(1-t^{\theta (j-i)}x_{i}y_{j}) 
\end{eqnarray}
In the second line we have used the properties (\ref{Us=-U})
and (\ref{Ux}).

We need to determine $[ x^{\sigma -1}]\,
\prod_{i \neq j}(1-t^{\theta (j-i)}x_{i}y_{j}) $
 up to equivalence under 
$\, {\acute{\sim}}_{y} \,$, which means finding the coefficient
of $x^{\sigma -1}$ neglecting any $y^{\eta}$ terms with
repeated parts. This coefficient must be a linear combination
of $y^{\sigma' -1}$ with $\sigma' \in S_{n}$.
Write
\begin{eqnarray} \nonumber
\prod_{i \neq j}(1-t^{\theta (j-i)}x_{i}y_{j})
&=&(1-x_{\sigma^{-1}(n)}y_{1})\ldots
(1-x_{\sigma^{-1}(n)}y_{\sigma^{-1}(n)-1})
(1-tx_{\sigma^{-1}(n)}y_{\sigma^{-1}(n)+1})
\ldots  
\\ &&\qquad\qquad 
\ldots (1-tx_{\sigma^{-1}(n)}y_{n})
\prod_{i \neq j,\,\sigma^{-1}(n)}(1-t^{\theta (j-i)}x_{i}y_{j})
\end{eqnarray}
It is clear that the coefficient of 
$x^{\sigma -1}=x_{\sigma^{-1}(n)}^{n-1}\ldots x_{\sigma^{-1}(1)}^{0} $ 
must be a linear combination
of $y^{\eta}$ with $\eta_{j} \geq 1$ for all $j$ except 
$j=\sigma^{-1}(n)$. Hence the coefficient of $x^{\sigma -1}$
must be a linear combination
of $y^{\eta}$ with
$\eta_{\sigma^{-1}(n)}=0$.
Continuing in this vein we reach the conclusion that the
coefficient of $x^{\sigma -1}$
 up to equivalence under 
$\, {\acute{\sim}}_{y} \,$
is a scalar multiple of $y^{\underline{\sigma}-1}$.
Hence
\begin{eqnarray}  \nonumber
[x^{\sigma -1}]\,\prod_{i \neq j}(1-t^{\theta (j-i)}x_{i}y_{j})
& {\acute{\sim}}_{y} &\left(
[x^{\sigma -1}y^{\underline{\sigma} -1}]\,
\prod_{i \neq j}(1-t^{\theta (j-i)}x_{i}y_{j})
\right)y^{\underline{\sigma} -1}
\\  \nonumber
&=&(-1)^{n(n-1)/2}t^{n-\sigma^{-1}(n)}\times
t^{n-\sigma^{-1}(n-1)-\theta(\sigma^{-1}(n)-\sigma^{-1}(n-1))}\ldots
\\  \nonumber
& & \qquad\quad \ldots\times 
t^{n-\sigma^{-1}(1)-\sum_{i=1}^{n}
\theta\left(\sigma^{-1}(i)-\sigma^{-1}(1)\right)}
y^{\underline{\sigma} -1}
\\   \nonumber
&=&(-t)^{n(n-1)/2}t^{-\#\{(i,j):\sigma^{-1}(j)>\sigma^{-1}(i),i<j\} }
y^{\underline{\sigma} -1}
\\  \label{Cauchy5}
&=&(-t)^{n(n-1)/2}t^{-l(\underline{\sigma})}
y^{\underline{\sigma} -1}
\end{eqnarray}
The stated result now follows after substituting (\ref{Cauchy5})
into (\ref{Cauchy4}), making use of (\ref{Cauchy1.5}),
and substituting the resulting identity in (\ref{Cauchy3}).
\hfill $\Box$ \medskip

We can now give a new derivation of the symmetric analogue of
Proposition \ref{E-sum}.
\begin{prop} {\rm\cite{Ma95}}\label{P-sum}
 We have
\begin{equation}  \label{P-sum1}
\Pi(x,y;q,t)
=
\sum_{\ka}\frac{1 }{ v_{\ka}(q,t)}P_{\ka}(x;q,t)\,P_{\ka}(y;q,t),
\quad \Pi(x,y;q,t):=\prod_{i,j=1}^{n}\frac{1 }{ (x_{i}y_{j};q)_{\lambda}}
\end{equation}
for scalars $v_{\ka}(q,t)$ independent of $n$.
\end{prop}
\noindent{\it Proof.} \quad
To derive (\ref{P-sum1}) we apply ${U^{-}}^{(x)}$
followed by ${U^{-}}^{(y)}|_{t\rightarrow t^{-1}}$
to both sides of (\ref{E-sum1}).
Write
\begin{equation}
\Omega(x,y;q,t)=\prod_{i,j=1}^{n}\frac{1 }{ (x_{i}y_{j};q)_{\lambda}}
 \prod_{i=1}^{n}\frac{1}{ 1-tx_{i}y_{i}}
\prod_{j<i}^{n}\frac{ 1-x_{i}y_{j}}{  1-tx_{i}y_{j}}
\end{equation}
Applying ${U^{-}}^{(x)}$ to
the left hand side of (\ref{E-sum1})
and using Lemma \ref{Cauchy} 
 we get
\begin{equation}  \label{P-sum2}
{U^{-}}^{(x)}\Omega(x,y;q,t)
=F(y)\,\Delta_{t}(x)\Pi(x,y;q,qt)
\end{equation}
Noting that $U^{-}U^{-}=[n]_{t^{-1}}!\,U^{-}$ and using
(\ref{Ux}) we have
\begin{equation}
{U^{-}}^{(y)}|_{t\rightarrow t^{-1}}\,F(y)
={U^{-}}^{(y)}|_{t\rightarrow t^{-1}}\,\Delta_{t^{-1}}(y)
=[n]_{t}!\,\Delta_{t^{-1}}(y)
\end{equation}
So applying ${U^{-}}^{(y)}|_{t\rightarrow t^{-1}}$ to
(\ref{P-sum2}) we obtain
\begin{equation}  \label{P-sum3}
{U^{-}}^{(y)}|_{t\rightarrow t^{-1}}\,{U^{-}}^{(x)}\,\Omega(x,y;q,t)
=
[n]_{t}!\,\Delta_{t}(x)\Delta_{t^{-1}}(y)
\Pi(x,y;q,qt)
\end{equation}
Applying 
${U^{-}}^{(y)}|_{t\rightarrow t^{-1}}{U^{-}}^{(x)}$ now
to the right hand side of (\ref{E-sum1}) and using (\ref{Us})
and Lemma \ref{s-qa} gives
\begin{eqnarray} \nonumber
&&\sum_{\rho}\frac{1}{ u_{\rho}(q,t)}{U^{-}}^{(x)}E_{\rho}(x;q,t)\,
{U^{-}}^{(y)}|_{t\rightarrow t^{-1}}
E_{\rho}(y;q^{-1},t^{-1})  \qquad\qquad\qquad\qquad
\\
&&\qquad\qquad\qquad \label{P-sum4}
={\sum_{\rho}}^{*}\frac{\beta_{\rho}(q,t)\beta_{\rho}(q^{-1},t^{-1})
}{ u_{\rho}(q,t)}
\Delta_{t}(x)\,P_{\rho^{+}-\delta}(x;q,qt)\,\Delta_{t^{-1}}(y)
P_{\rho^{+}-\delta}(x;\frac{1}{ q},\frac{1}{ qt})
\end{eqnarray}
where the $*$ denotes that the sum is restricted to
$\rho$ with distinct parts.
Equating (\ref{P-sum3}) and (\ref{P-sum4}), using (\ref{4.1.1})
and letting $qt\mapsto t$ we obtain (\ref{P-sum1}).
The stability property of the symmetric Macdonald polynomials
\cite{Ma95}
\begin{equation}
P_{\ka}(x_1,\ldots,x_{n-1},0;q,t)=
\left\{
\begin{array}{ll} 
 P_{\ka}(x_1,\ldots,x_{n-1};q,t)  & {\rm length}(\ka)\leq n-1
\\0 & \mbox{else}
\end{array}
\right.
\end{equation}
applied to (\ref{P-sum1}) shows that the
$v_{\ka}(q,t)$ are independent of $n$.
\hfill $\Box$ \medskip

Define the polynomials $g_{\ka}(x;q,t)$ by \cite{Ma95}
\begin{equation}
\Pi(x,y;q,t):=\sum_{\ka}g_{\ka}(x;q,t)m_{\ka}(y)
\end{equation}

\begin{cor} {\rm cf. \cite[p310-11,313]{Ma95}}\quad
Define an scalar product by 
$\langle P_{\ka}(x;q,t),P_{\mu}(x;q,t)\rangle_{g}
:=\nu_{\ka}(q,t)\delta_{\mu\ka}$. We have
\begin{equation}
\langle g_{\mu}(x;q,t),m_{\ka}(x)\rangle_{g}=\delta_{\mu\ka}
\end{equation}
and hence the $g_{\mu}(x;q,t)$ are a basis for the multivariable
polynomials with coefficients in ${\mathbb{Q}}(q,t)$. 
\end{cor}
\noindent {\it Proof.}
Similar to the proof of Corollary \ref{<>q1}.
\hfill $\Box$ \medskip

In order to proceed further with the development of
the symmetric theory we require the following
symmetrization formulas.
\begin{lemma} \label{coef}
Let $\eta^{R}:= \underline{(\eta^{+})}$ and 
$f_{j}:=f_{j}(\eta):=\#\{ i:\eta_{i}=j\}$. Then
\begin{eqnarray} \label{coef.1}
\mbox{a)}&\qquad& P_{\eta^{+}}(x;q,t)=
t^{-n(n-1)/2}\prod_{j=0}^{\eta_{1}^{+}}\frac{1}{ [f_{j}]_{t^{-1}}!}
U^{+}E_{\eta^{R}}(x;q,t)
\\
\mbox{b)} &\qquad& S_{\eta^{+}}(x;q,t)=(-1)^{n(n-1)/2}
U^{-}E_{\eta^{R}}(x;q,t)
\end{eqnarray}
\end{lemma}
\noindent {\it Remark.}\quad
Using the theory of the symmetric Macdonald polynomials
Baker and Forrester \cite[(5.8),(5.18)]{BF97g}   have
derived a more general formula for the constant relating
$U^{+}E_{\eta}$ and $P_{\eta^{+}}$. Their expression is
not in the same form as (\ref{coef.1}), although they can be shown to
be equal using the first equality of (\ref{P(t).0}).
\medskip

\noindent {\it Proof.}
We shall only consider (a) as the proof of (b) is similar.
From the triangular structure of $E_{\eta^{R}}(x;q,t)$
\begin{equation}
U^{+}E_{\eta^{R}}(x;q,t)=U^{+}x^{\eta^{R}}+\sum_{\nu<\eta^{+}}
a_{\nu}m_{\nu}(x)
\end{equation}
From (\ref{U+E}) we know that $U^{+}E_{\eta^{R}}$
is a scalar multiple of $P_{\eta^{+}}$. To find the 
scalar multiple we need to determine
$[m_{\eta^{+}}]\,U^{+}x^{\eta^{R}}$.
Suppose $\sigma =s_{i_{1}}\ldots s_{i_{p}}$ is the reduced
decomposition of the permutation $\sigma$. For all $k=1,\ldots,p$
\begin{equation}
(s_{i_{k+1}}\ldots s_{i_{p}}\eta^{R})_{k}\leq
(s_{i_{k+1}}\ldots s_{i_{p}}\eta^{R})_{k+1}
\end{equation}
From the action of $T_i$ (\ref{TE}) we then have
\begin{equation}
T_{\sigma}x^{\eta^{R}}=t^{l(\sigma)}x^{\sigma\eta^{R}}
+\sum_{\mu\prec\sigma\eta^{R}}b_{\mu}x^{\mu}
\end{equation}
It follows that
\begin{eqnarray} \nonumber
[x^{\eta^{+}}]\,U^{+}x^{\eta^{R}}&=&
[x^{\eta^{+}}]\sum_{\sigma\eta^{R}=\eta^{+}}T_{\sigma}x^{\eta^{R}}
= \sum_{\sigma\eta^{R}=\eta^{+}}t^{l(\sigma)}
= \sum_{\sigma\eta^{R}=\eta^{R}}t^{l(\underline{\sigma})}
=t^{n(n-1)/2}\sum_{\sigma\eta^{R}=\eta^{R}}t^{-l(\sigma)}
\\ 
&=&  \label{coef.2}
t^{n(n-1)/2}\prod_{j=0}^{\eta^{+}_{1}}
\sum_{\sigma\in S_{f_{j}}}t^{-l(\sigma)}
= t^{n(n-1)/2}\prod_{j=0}^{\eta^{+}_{1}}[f_{j}]_{t^{-1}}!
\end{eqnarray}
Since $U^{+}x^{\eta^{R}}$ is symmetric this shows
that the coefficient of $m_{\eta^+}$ in $U^{+}x^{\eta^{R}}$
is given by the right hand side of (\ref{coef.2})
as required by (\ref{coef.1}).
\hfill $\Box$ \medskip

We can now deduce the symmetric analogue of Proposition \ref{E(t)}.
\begin{prop} {\rm\cite{Ma95}} \label{P(t)}
\begin{equation}  \label{P(t).0}
P_{\eta^{+}}(t^{\underline{\delta}};q,t)
=\frac{t^{l(\eta^{R})} \, [n]_{t^{-1}}! }{ 
\prod_{i=0}^{\eta_{1}^{+}}  \,[f_{j}]_{t^{-1}}!       }
\frac{e_{\eta^{R}}(q,t)}{ d_{\eta^{R}}(q,t)}
=t^{l(\eta^+)}\frac{b_{\eta^{+}}(q,t)}{ h_{\eta^{+}}(q,t)}
\end{equation}
where
\begin{equation}
 h_{\ka}(q,t):=\prod_{s\in\ka}(1-q^{a(s)}t^{l(s)+1})
\end{equation}
\end{prop}
\noindent {\it Proof.}  \quad
Applying Lemma \ref{coef}(a) we have
\begin{equation}
P_{\eta^{+}}(t^{\underline{\delta}};q,t)
=\frac{t^{-n(n-1)/2}}{\prod_{j=0}^{\eta_{1}^{+}} [f_{j}]_{t^{-1}}!}
\sum_{\sigma\in S_{n}}
\left. \left( T_{\sigma} E_{\eta^{R}}(x;q,t)
\right) \right|_{x=t^{\underline{\delta}}}
\end{equation}
Using (\ref{Tf}) we obtain
\begin{equation}
P_{\eta^{+}}(t^{\underline{\delta}};q,t)
=\frac{t^{-n(n-1)/2}\sum_{\sigma\in S_{n}}t^{l(\sigma)}
}{\prod_{j=0}^{\eta_{1}^{+}}
[f_{j}]_{t^{-1}}!}E_{\eta^{R}}(t^{\underline{\delta}};q,t)
\end{equation}
Since $t^{-n(n-1)/2}\sum_{\sigma\in S_{n}}t^{l(\sigma)}
=\sum_{\sigma\in S_{n}}t^{l(\underline{\sigma})}
=[n]_{t^{-1}}!$
we obtain the first equality in (\ref{P(t).0})
by using Proposition \ref{propE(t)}.

The second equality follows immediately from the identities
\begin{equation}  \label{P(t).1}
\frac{[n]_{t}!}{\prod_{j=0}^{\eta_{1}^{+}}
[f_{j}]_{t}!}
\frac{e_{\eta^{R}}(q,t)}{ d_{\eta^{R}}(q,t)}
=\frac{b_{\eta^{+}}(q,t)}{ h_{\eta^{+}}(q,t)}
\end{equation}
\begin{equation}   \label{P(t).2}
t^{l(\eta^{R})-l(\eta^{+})}
=\frac{[n]_{t}!}{[n]_{t^{-1}}!}
\prod_{j=0}^{\eta_{1}^{+}}
\frac{[f_{j}]_{t^{-1}}!}{[f_{j}]_{t}!}
\end{equation}
For the first identity we use (\ref{e-G}) and (\ref{gamma1})
to obtain 
\begin{eqnarray}
\frac{e_{\eta^{R}}(q,t)}{ b_{\eta^{+}}(q,t)}
&=&\frac{1 }{ [n]_{t}!}\prod_{i=1}^{n}
[{{\lambda}^{-1}\eta_{i}^{+} }+n-i+1]_{t}
\\
&=&\frac{[f_{0}(\eta)]_{t}! }{ [n]_{t}!}
\prod_{i=1}^{n-f_{0}(\eta)}[{{\lambda}^{-1}\eta_{i}^{+} }+n-i+1]_{t}
\end{eqnarray}
It suffices then to show that
\begin{equation} \label{P(t).3}
\frac{ h_{\eta^{+}}(q,t)}{ \prod_{j=0}^{\eta_{1}^{+}}
[f_{j}]_{t}!}
=\frac{d_{\eta^{R}}(q,t)}{ 
\prod_{i=1}^{n-f_{0}(\eta)}[{{\lambda}^{-1}\eta_{i}^{+} }+n-i+1]_{t}}
\end{equation}
This is an easy consequence of a natural $q$-generalisation of
the argument used in \cite{BF97c} to prove the corresponding identity
in the Jack polynomial theory.

We now turn to the second identity. Noting that 
$\frac{[m]_{t}!}{ [m]_{t^{-1}}!} =t^{m(m-1)/2}$ we have
\begin{equation}  
\frac{[n]_{t}!}{[n]_{t^{-1}}!}
\prod_{j=0}^{\eta_{1}^{+}}
\frac{[f_{j}]_{t^{-1}}!}{[f_{j}]_{t}!}
=t^{n(n-1)/2-\sum_{j=0}^{\eta_{1}^{+}}f_{j}(f_{j}-1)/2}
\end{equation}
It follows from Lemma \ref{lem2.1} that 
\begin{equation}
l(\eta^{R})=l(\sigma)+l(\eta^{+})
\end{equation}
 where $\sigma$ is the permutation of minimum length for which
$\eta^{+}=\sigma(\eta^{R})$. Since the minimum such length is
$l(\sigma)
=n(n-1)/2-\sum_{j=0}^{\eta_{1}^{+}}f_{j}(f_{j}-1)/2$
we obtain (\ref{P(t).2}).
\hfill $\Box$ \medskip

\begin{prop}{\rm\cite{Ma95}}
Let ${\cal N}_{\ka}^{(P)}(q,t):=\langle P_{\ka}(x;q,t),P_{\ka}(x;q,t)
\rangle_{q,t}$. With $\eta^+ =\ka$ we have
\begin{equation}  \label{NP.0}
\frac{{\cal N}_{\eta^+}^{(P)}(q,t)}{ {\cal N}_{0}^{(P)}(q,t)}
=\frac{[n]_{t}!}{\prod_{j=0}^{\eta_{1}^{+}}
[f_{j}]_{t}!}
\frac{{d'}_{\eta^+}(q,t)\,e_{\eta^{R}}(q,t)}{
d_{\eta^{R}}(q,t)\,{e'}_{\eta^R}(q,t)}
=\frac{b_{\eta^{+}}(q,t)\, {d'}_{\eta^+}(q,t)
}{ h_{\eta^{+}}(q,t){e'}_{\eta^+}(q,t)}
\end{equation}
\end{prop}
\noindent {\it Proof.}  \quad
We have
\begin{eqnarray}  \nonumber
\langle U^{+}E_{\eta^{R}},U^{+}E_{\eta^{R}}\rangle_{q,t}
&=&\sum_{\sigma \in S_{n}}
\langle U^{+}E_{\eta^{R}},T_{\sigma}E_{\eta^{R}}\rangle_{q,t}
=
\sum_{\sigma \in S_{n}}
\langle T_{\sigma}^{-1}U^{+}E_{\eta^{R}},E_{\eta^{R}}\rangle_{q,t}
\\  \label{NP.1}
&=&\sum_{\sigma \in S_{n}}t^{-l(\sigma)}
\langle U^{+}E_{\eta^{R}},E_{\eta^{R}}\rangle_{q,t}
=
[n]_{t^{-1}}\,
\langle U^{+}E_{\eta^{R}},E_{\eta^{R}}\rangle_{q,t}
\end{eqnarray}
In the second equality we have used the fact that $T_{i}^{-1}$ is the 
adjoint operator of $T_{i}$ while in the third equality we have used
(\ref{TU+}). Multiplying each side of (\ref{NP.1}) by
$\prod_{i=0}^{\eta_{1}^{+}}1/[f_{j}]_{t}![f_{j}]_{t^{-1}}!$ and using
Proposition \ref{coef} we obtain
\begin{equation}
\langle P_{\eta^+},P_{\eta^+}\rangle_{q,t}
=\frac{t^{n(n-1)/2}\,[n]_{t^{-1}}!}{
\prod_{i=0}^{\eta_{1}^{+}}  \,[f_{j}]_{t}!  }
\langle P_{\eta^+},E_{\eta^R}\rangle_{q,t}
\end{equation}
Using (\ref{P-E}) and the orthogonality of the 
non-symmetric Macdonald polynomials we get
\begin{equation}
\langle P_{\eta^+},P_{\eta^+}\rangle_{q,t}
\frac{[n]_{t}! }{\prod_{i=0}^{\eta_{1}^{+}}  \,[f_{j}]_{t}!}
\frac{ {d'}_{\eta^+}(q,t) }{ {d'}_{\eta^R}(q,t)}
\langle E_{\eta^R},E_{\eta^R}\rangle_{q,t}
\end{equation}
Dividing each side by ${\cal N}_{0}^{(P)}(q,t)={\cal N}_{0}^{(E)}(q,t)$
and using Proposition \ref{NE} we obtain the equality on
the right hand
side of (\ref{NP.0}). The second identity follows from
using the identity (\ref{P(t).1}).
\hfill $\Box$ \medskip

It remains to establish the analogue of Proposition \ref{Bi-E}
and to specify the constant $v_{\ka}(q,t)$ appearing in Proposition
\ref{P-sum}. We proceed as in the derivation of Proposition \ref{Bi-E}
using (\ref{P-sum1}), (\ref{P(t).0}) and the identity
\begin{equation}
\frac{ b_{\eta^+}(\frac{1}{ q},\frac{1}{ t}) 
}{ 
h_{\eta^+}(\frac{1}{ q},\frac{1}{ t})}
=t^{\left( l(\eta^+ ) +l'(\eta^+)-(n-1)|\eta| \right)}
\frac{ b_{\eta^+}(q,t) }{ h_{\eta^+}(q,t)}
\end{equation}
We obtain 
\begin{equation}  \label{Bi-P.-1}
\prod_{i=1}^{n}\frac{1}{ (x_{i};q)_{r}}=
\sum_{\eta^+}
\frac{[q^{r}]_{\eta^+} }{ v_{\eta^+}(q,t) h_{\eta^+}(q,t)}
P_{\eta^+} (x;q,t)
\end{equation}
Now substituting (\ref{P-E}) for $P_{\eta^+}$ and comparing the 
results with (\ref{Bi-E.1}) we can read off the value of
$v_{\eta^+}(q,t)$.
\begin{prop} {\rm\cite{MN}}
\begin{equation}
v_{\ka}(q,t) =\frac{{d'}_{\ka}(q,t)}{ h_{\ka}(q,t) }
\end{equation}
\end{prop}
Substituting this result back into (\ref{Bi-P.-1}) we obtain
the $q$-binomial theorem involving the symmetric Macdonald polynomials.
\begin{prop}{\rm\cite{Kan96}}
\begin{equation}  \label{Bi-P}
\prod_{i=1}^{n}\frac{1}{ (x_{i};q)_{r}}=
\sum_{\ka}
\frac{[q^{r}]_{\ka} }{  {d'}_{\ka}(q,t)}
P_{\ka} (x;q,t)
\end{equation}
\end{prop}

To tie things up we shall prove that 
$F(y)$ appearing in (\ref{Cauchy1}) is equal to
$\Delta_{t^{-1}}(y)$ and hence derive a 
$q$-generalisation of the Cauchy double
alternant formula. The derivation will
also yield the value of the constant
$\beta_{\eta}(q,t)$ appearing in (\ref{Us}).
\begin{lemma} \label{2.Cauchy}
\begin{eqnarray}  \label{2.Cauchy.a}
\mbox{a)} &\qquad &
F(y)=\Delta_{t^{-1}}(y)
\\
\mbox{b)} \label{2.Cauchy.b}&\qquad &
\beta_{\sigma (\eta^+)}(q,t)= (-1)^{l(\sigma)}
\frac{ {d'}_{\sigma (\eta^+)}(q,t) \,{h}_{(\eta^+ -\delta)}\,(q,qt)
}{ {d}_{\eta^+}(q,t) \,{d'}_{(\eta^+ -\delta)}\,(q,qt)}
\end{eqnarray}
\end{lemma}

\noindent {\it Proof.}
As in the proof of Proposition \ref{P-sum} apply
${U^{-}}^{+}$ to (\ref{E-sum1}) and cancel the
factor $\Delta_{t}(x)$ from the resulting expression.
Substituting the right hand side of (\ref{E-sum1})
with $t\rightarrow qt$, 
using (\ref{4.1.1}) and multiplying by $\Delta_{t^{-1}}(y)$,
we obtain
\begin{equation}  \label{2.Cauchy.1}
F(y)\sum_{\ka}\frac{1}{ v_{\ka}(q,qt)}P_{\ka}(x;q,qt)
\Delta_{t^{-1}}(y) P_{\ka}(y;\frac{1}{ q},\frac{1}{ qt})
=\Delta_{t^{-1}}(y) {\sum_{\eta}}^{*}
\frac{t^{n(n-1)/2} \beta_{\eta}(q,t)}{ u_{\eta}(q,t)}
P_{(\eta^{+}-\delta)}(x;q,qt)E_{\eta}(y;\frac{1}{ q},\frac{1}{ t})
\end{equation}
Using (\ref{s-qa1}) and (\ref{P-E}) we can write
\begin{equation}
\Delta_{t^{-1}}(y)P_{\ka}(y;\frac{1}{ q},\frac{1}{ qt})
=t^{n(n-1)/2}\sum_{\sigma\in S_{n}}(-t)^{l(\sigma)}
\frac{d_{\sigma (\eta^+)}(\frac{1}{ q},\frac{1}{ t})
}{
d_{ \eta^+}(\frac{1}{ q},\frac{1}{ t})}
E_{\sigma (\eta^+)}(y;\frac{1}{ q},\frac{1}{ t})
\end{equation}
Substituting into the left hand side of (\ref{2.Cauchy.1}) and
equating the coefficients of $ P_{\ka}(y;\frac{1}{ q},\frac{1}{ qt})$
we obtain for $\eta$ with no repeated parts
\begin{equation}  \label{2.Cauchy.2}
F(y)\sum_{\sigma\in S_{n}}
\frac{  (-t)^{l(\sigma)} 
d_{\sigma (\eta^+)}(\frac{1}{ q},\frac{1}{ t}) }{
v_{(\eta^{+}-\delta)}(q,qt) d_{ \eta^+}(\frac{1}{ q},\frac{1}{ t})}
E_{\sigma (\eta^+)}(y;\frac{1}{ q},\frac{1}{ t})
=
\Delta_{t^{-1}}(y) \sum_{\sigma\in S_{n}}
\frac{\beta_{\sigma (\eta^+)}(q,t)}{ u_{\sigma (\eta^+)}(q,t)}
E_{\sigma (\eta^+)}(y;\frac{1}{ q},\frac{1}{ t})
\end{equation}
Let $\bar{\,}$ be the linear operator such that
\begin{equation}
\overline{x^{\mu}}= \left\{ \begin{array}{ll}
x^{\mu} &\quad   \mbox{$\mu$ has no repeated parts}
\\
0 &\quad \mbox{else}
\end{array}
\right.
\end{equation}
Apply this operator to (\ref{2.Cauchy.2}). Since
$\overline{F}(y)={\overline{\Delta}}_{t^{-1}}(y)$,
equating the coefficients of the linearly independent
$\{{\overline{\Delta}}_{t^{-1}}(y)
{\overline{E}}_{\sigma (\eta^{+})}
(y;q^{-1},t^{-1})\}$ results in
\begin{equation} \label{2.Cauchy.3}
\beta_{\sigma (\eta^+)}(q,t)=
(-t)^{l(\sigma)}\frac{u_{\sigma (\eta^+)}(q,t) }{
v_{(\eta^{+}-\delta)}(q,qt)}
\frac{d_{\sigma (\eta^+)}(\frac{1}{ q},\frac{1}{ t}) 
}{
 d_{ \eta^+}(\frac{1}{ q},\frac{1}{ t})}
\end{equation}
Substituting back into (\ref{2.Cauchy.2}) we obtain
 (\ref{2.Cauchy.a}).
To obtain (\ref{2.Cauchy.b}) we simplify (\ref{2.Cauchy.3}) using 
the identity
\begin{equation}
\frac{ d_{ \eta^+}( q, t)\, d_{\sigma (\eta^+)}(\frac{1}{ q},\frac{1}{ t})
}{ d_{\sigma (\eta^+)}( q, t)
d_{\eta^+}(\frac{1}{ q},\frac{1}{ t})}
=t^{-l(\sigma )}
\end{equation}
\hfill $\Box$ \medskip

\noindent {\it Remark.}\quad The expression for the constant
$\beta_{\sigma (\eta^+)}(q,t)$ in Lemma \ref{2.Cauchy} can
be simplified using a natural $q$-generalisation of the
argument in \cite{BF97c}. 
The simplification gives
\begin{equation}
\beta_{\sigma (\eta^+)}(q,t)= (-1)^{l(\sigma)}
\frac{{d'}_{\sigma (\eta^{+})}(q,t) }{ 
{d'}_{\sigma (\eta^{R})}(q,t)}
\end{equation}

\end{document}